\documentclass{article}
\usepackage{latexsym}
\usepackage{epsfig}
\usepackage{amsmath}
\usepackage{amssymb}
\usepackage{amsthm}
\newtheorem{theorem}{Theorem}[section]
\newtheorem{corollary}[theorem]{Corollary}

\newtheorem{lemma}[theorem]{Lemma}

\theoremstyle{remark}

\theoremstyle{definition}

\DeclareMathOperator\diag{diag}

\begin{document}

\title{Optimal inequalities between distances \\ in convex projective domains}

\author{Roland Hildebrand \thanks{%
Univ.\ Grenoble Alpes, CNRS, Grenoble INP, LJK, 38000 Grenoble, France
({\tt roland.hildebrand@univ-grenoble-alpes.fr}).}}

\maketitle

\begin{abstract}
On any proper convex domain in real projective space there exists a natural Riemannian metric, the Blaschke metric. On the other hand, distances between points can be measured in the Hilbert metric. Using techniques of optimal control, we provide inequalities lower bounding the Riemannian length of the line segment joining two points of the domain by the Hilbert distance between these points, thus strengthening a result of Tholozan. Our estimates are valid for a whole class of Riemannian metrics on convex projective domains, namely those induced by convex non-degenerate centro-affine hypersurface immersions. If the immersions are asymptotic to the boundary of the convex cone over the domain, then we can also upper bound the Riemmanian length. On these classes, and in particular for the Blaschke metric, our inequalities are optimal.
\end{abstract}

Keywords: Hilbert distance, Blaschke metric, centro-affine hypersurface immersion, affine hypersphere

MSC 2020: 53A15, 52A38, 58E10

\section{Introduction}

Proper convex domains $\Omega \subset \mathbb RP^n$ allow to define several projectively invariant distances. On of them is the Hilbert distance. For two distinct points $a,b \in \Omega$, it is defined by
\[ d^H(a,b) = \frac12\log\frac{||ya||\cdot||xb||}{||yb||\cdot||xa||} = \frac12\log(a,b;y,x),
\]
where $x,y \in \partial\Omega$ are the boundary points lying on the projective line $l$ through $a,b$, such that the order of points on $l$ is $x,a,b,y$, and the quantities $||ya||,\dots$ are the corresponding coordinate differences in any affine chart on $l$. The quantity $(a,b;x,y)$ is the projective cross-ratio of the quadruple of points.

Another possibility to define distances in the domain $\Omega$ is by a Riemannian metric. Let $K_{\Omega} \subset \mathbb R^{n+1}$ be the proper convex cone over the closure of the domain $\Omega$. Lift $\Omega$ into the interior of $K_{\Omega}$. Such a hypersurface immersion $f: \Omega \to K_{\Omega}^o$, smooth enough and equipped with the position vector as transversal vector field, gives rise to a centro-affine fundamental form $h$ on $\Omega$. If the immersion is locally strongly convex, then this form defines a Riemannian metric on $\Omega$. We suppose throughout the paper that the immersion is of hyperbolic type, i.e., the hypersurface is bent away from the origin.

As a special case, we obtain the \emph{Blaschke metric} (or Cheng-Yau metric) if the immersion is a proper affine hypersphere with mean curvature $-1$ which is asymptotic to the boundary $\partial K_{\Omega}$. The Blaschke metric gives rise to the \emph{Blaschke distance} $d^B$ on $\Omega$ (see \cite{Loftin01},\cite{BenoistHulin13},\cite{KimPapadopoulos14},\cite{LoftinMacintosh09}).

Along with the fundamental form, a centro-affine hypersurface immersion $f$ defines a symmetric 3-rd order tensor field on $\Omega$, the \emph{cubic form} $C$. On affine hyperspheres the cubic form is bounded by an explicit function of the dimension $n$ \cite{Hildebrand13a},
\begin{equation} \label{C_bound_AHS}
C(X,X,X) \leq 2\frac{n-1}{\sqrt{n}}(h(X,X))^{3/2}
\end{equation}
for all tangent vector fields $X$ on $\Omega$.

In \cite{BenoistHulin13}, Benoist and Hulin proved by a general compactness argument \cite{Benzecri60} that both distances $d^B,d^H$ are strongly equivalent, i.e., one can be bounded by a multiple of the other, where the constants depend only on the dimension $n$. In \cite{Tholozan17}, Tholozan used this result to prove the remarkable inequality
\[ d^B(a,b) < d^H(a,b) + 1
\]
for all pairs of points $a,b \in \Omega$. This inequality actually holds also for Riemannian distances $d^R$ generated by general convex non-degenerate centro-affine hypersurface immersions of $\Omega$.

In this paper we improve this result to optimality as follows.

{\theorem \label{thm:C_libre} Let $\Omega \subset \mathbb RP^n$ be a proper convex domain, and let $f: \Omega \to K_{\Omega}^o \subset \mathbb R^{n+1}$ be a non-degenerate convex lift of class $C^2$ into the interior of the convex cone over $\Omega$. Let $h$ be the centro-affine metric induced on $\Omega$ by $f$, and let $d^R$ be the corresponding geodesic distance. Let further $d^H$ be the Hilbert distance on $\Omega$. Then for every pair of points $a,b \in \Omega$ the inequalities
\[ d^R(a,b) \leq l^R(a,b) < \log\left( \exp(d^H(a,b)) + \sqrt{\exp(2d^H(a,b)) - 1} \right) < d^H(a,b) + \log 2
\]
hold, where $l^R(a,b)$ is the Riemannian length of the projective line segment joining the points $a,b$. The inequalities cannot be improved. }

\medskip

In the case when $\Omega$ is an ellipsoid both the Hilbert and the Blaschke metric coincide, and $\Omega$ is isometric to hyperbolic space \cite{Tholozan17}. Note that in this case the cubic form of the affine sphere over $\Omega$ vanishes. This suggests that the deviation of the Hilbert distance from the centro-affine Riemannian distance can somehow be controlled by a bound on the cubic form of the centro-affine hypersurface immersion $f$. We investigate this dependence and provide optimal inequalities between $l^R(a,b)$ and $d^H(a,b)$ with the bound on the cubic form appearing as a parameter.

{\theorem \label{thm:C_bound} Let $\Omega \subset \mathbb RP^n$ be a proper convex domain, and let $f: \Omega \to K_{\Omega}^o \subset \mathbb R^{n+1}$ be a non-degenerate convex lift of class $C^3$ into the interior of the convex cone over $\Omega$. Let $C$ be the cubic form and $h$ the centro-affine metric induced on $\Omega$ by $f$, and let $d^R$ be the corresponding geodesic distance. Suppose that the cubic form satisfies the bound
\begin{equation} \label{Cbound}
C(X,X,X) \leq 2\gamma(h(X,X))^{3/2}
\end{equation}
for all tangent vector fields $X$ on $\Omega$, where $\gamma > 0$ is some constant. Let further $d^H$ be the Hilbert distance on $\Omega$. Then for every pair of points $a,b \in \Omega$ the inequalities
\begin{align*}
d^R(a,b) &\leq l^R(a,b) < \left\{ \begin{array}{rcl} \frac{2\mu}{\mu^2+1}\log\frac{(\mu^2+1)(E - 1) + 2}{2},&\quad& d^H(a,b) \leq \log\frac{\mu^2+1}{\mu^2-1}, \\
\log\frac{(\mu-1)(\sqrt{E+1} + \sqrt{E-1})}{(\mu+1)(\sqrt{E+1} - \sqrt{E-1})} + \frac{2\mu}{\mu^2+1}\log\frac{2\mu^2}{\mu^2-1},&& d^H(a,b) \geq \log\frac{\mu^2+1}{\mu^2-1}, \end{array} \right. \\
&< d^H(a,b) + \log 2 - \log\frac{\mu + 1}{\mu - 1} + \frac{2\mu}{\mu^2+1}\log\frac{2\mu^2}{\mu^2-1}
\end{align*}
hold, where $l^R(a,b)$ is the Riemannian length of the projective line segment joining the points $a,b$, $E = \exp(d^H(a,b))$, and $\mu = \frac{\gamma}{2} + \sqrt{1 + \frac{\gamma^2}{4}}$. The inequalities cannot be improved. }

\medskip

Under the additional assumption that the centro-affine immersion is asymptotic to the boundary of the convex cone $K_{\Omega}$ over the domain we can derive also a lower bound on the Riemannian length.

{\theorem \label{thm:C_lower_bound} Assume the notations and conditions of Theorem \ref{thm:C_bound}, and suppose in addition that the immersion $f$ is asymptotic to the boundary $\partial K_{\Omega}$. Then the inequalities
\[ l^R(a,b) > \frac{2\mu}{\mu^2 + 1}\log\frac{\exp(d^H(a,b))(\mu^2 + 1) + \mu^2 - 1}{2\mu^2} > \frac{2\mu}{\mu^2 + 1}\left( d^H(a,b) - \log\frac{2\mu^2}{\mu^2 + 1} \right)
\]
hold for every $a,b \in \Omega$. The inequalities cannot be improved. }

\medskip

\begin{figure}
\centering
\includegraphics[width=16.17cm,height=6.86cm]{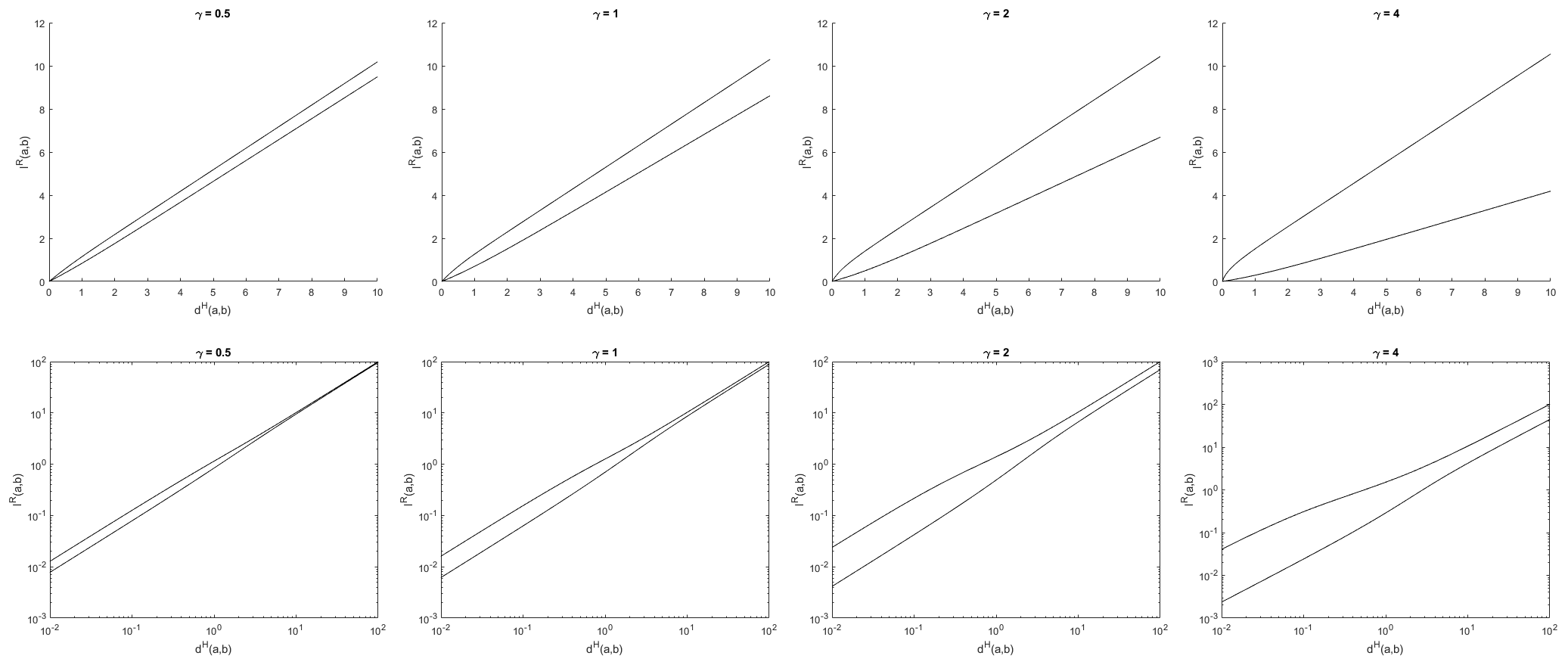}
\caption{Upper and lower bounds on the Riemannian length of line segments as a function of the Hilbert distance between the end-points for different values of $\gamma$. Linear plot (top two rows) and log-log plot (bottom two rows).}
\label{fig:bounds}
\end{figure}

The upper and lower bounds on $l^R(a,b)$ are depicted in Fig.~\ref{fig:bounds}. It is also possible to obtain bounds on the geodesic distance $d^R$.

{\theorem \label{thm:C_lower_bound_geodesic} Assume the notations and conditions of Theorem \ref{thm:C_lower_bound}. Then the inequalities
\[ \mu^{-1}d^H(a,b) \leq d^R(a,b) \leq \mu d^H(a,b)
\]
hold for every $a,b \in \Omega$. }

\medskip

In the case of affine spheres which are asymptotic to $\partial K_{\Omega}$ we get $\mu = \sqrt{n}$ by virtue of \eqref{C_bound_AHS}. Applying Theorems \ref{thm:C_bound}, \ref{thm:C_lower_bound}, \ref{thm:C_lower_bound_geodesic} yields the following result.

{\corollary Let $\Omega \subset \mathbb RP^n$ be a proper convex domain. Let $d^H(a,b)$ be the Hilbert distance, $d^B(a,b)$ the geodesic distance in the Blaschke metric, and $l^B(a,b)$ the Riemannian length of the straight line segment in the Blaschke metric between points $a,b \in \Omega$. Then
\begin{align*}
d^B(a,b) &\leq l^B(a,b) < \left\{ \begin{array}{rcl} \frac{2\sqrt{n}}{n+1}\log\frac{(n+1)(E - 1) + 2}{2},&\quad& d^H(a,b) \leq \log\frac{n+1}{n-1}, \\
\log\frac{(\sqrt{n}-1)(\sqrt{E+1} + \sqrt{E-1})}{(\sqrt{n}+1)(\sqrt{E+1} - \sqrt{E-1})} + \frac{2\sqrt{n}}{n+1}\log\frac{2n}{n-1},&& d^H(a,b) \geq \log\frac{n+1}{n-1}, \end{array} \right. \\
&< d^H(a,b) + \log 2 - \log\frac{\sqrt{n} + 1}{\sqrt{n} - 1} + \frac{2\sqrt{n}}{n+1}\log\frac{2n}{n-1},\\
l^B(a,b) &> \frac{2\sqrt{n}}{n + 1}\log\frac{E(n + 1) + n - 1}{2n} > \frac{2\sqrt{n}}{n + 1}\left( d^H(a,b) - \log\frac{2n}{n + 1} \right),
\end{align*}
\[ \frac{1}{\sqrt{n}}d^H(a,b) \leq d^B(a,b) \leq \sqrt{n}d^H(a,b),
\]
where $E = \exp(d^H(a,b))$. \qed }

This yields an explicit estimate of the constants realizing the equivalence of the two metrics.

\medskip

Connections between the affine metric of a complete hyperbolic affine sphere and the Hilbert metric have also been investigated in \cite{BenoistHulin14}. It has been established that one metric is Gromov hyperbolic if and only if the other is. Asymptotic properties of the Hilbert distance have been studied in \cite{Methou04}.

\medskip

The remainder of the paper is structured as follows. In Section \ref{sec:Bellman} we outline the proof strategy for the presented results. It is based on presenting an explicit solution of the Bellman equation. In Section \ref{sec:gen_ineq} we prove Theorem \ref{thm:C_libre}, in Section \ref{sec:bound_ineq} we prove Theorem \ref{thm:C_bound}, while in Section \ref{sec:lower_bound} we prove Theorems \ref{thm:C_lower_bound} and \ref{thm:C_lower_bound_geodesic}.

\section{Proof strategy: Bellman function} \label{sec:Bellman}

The bounds in Theorems \ref{thm:C_libre}, \ref{thm:C_bound} (Theorem \ref{thm:C_lower_bound}) result from maximizing (minimizing) a Riemannian length integral under some constraints, which can be cast in the form of an optimal control problem. In this section we describe a generic method to show the optimality of a solution to such a problem.

Suppose a controlled dynamical system $\dot x(t) = g(x(t),u(t))$ is given, where $x(t)$ is the sought scalar or vector-valued solution, and $u(t)$ is the control taking values in some set $U$, which may depend on $x$. Suppose further that $x$ is constrained to some closed convex set $X$. We seek to maximize an integral functional
\[ {\cal I}(x(\cdot)) = \int_{t_i}^{t_f} L(x(t),u(t))\,dt
\]
over the trajectories of the system, with initial and final conditions $(t_i,x(t_i)) \in M_i$, $(t_f,x(t_f)) \in M_f$, where $M_i,M_f$ are some subsets of $\mathbb R \times X$.

Such a problem may be solved by optimal control techniques, notably the Pontryagin Maximum Principle (PMP) \cite{PBGM62} or in the simple case of unconstrained dynamics by the Euler-Lagrange equation \cite{Gelfand63}. These tools provide necessary optimality conditions on the solution. They are thus efficient in finding potential solutions, but in order to actually prove optimality we have to use sufficient conditions. The simplest way to provably demonstrate the optimal value of the problem is to present the \emph{Bellman function} $B(t,x)$ \cite{Bellman}.

The value of the Bellman function at a point $(t_0,x_0) \in \mathbb R \times X$ is defined by the maximal value which can be achieved by the functional ${\cal I}$ over trajectories in $X$ with initial point $(t_0,x_0)$ and end-point $(t_f,x(t_f)) \in M_f$. In particular, the Bellman function is continuous and we have
\begin{equation} \label{Bellman_initial}
B(t,x) = 0\quad \forall\ (t,x) \in M_f.
\end{equation}
For every $x \in X$, let $U_x \subset U$ be the set of control values $u$ such that $g(x,u)$ is tangent to $X$ at $x$, i.e., application of a control satisfying $u \in U_x$ a.e.\ does not lead out of $X$. Then the Bellman function satisfies the \emph{Bellman equation}
\begin{equation} \label{Bellman_equation}
- \frac{\partial B}{\partial t} = \max_{u \in U_x} \left( \frac{\partial B}{\partial x}g(x,u) + L(x,u) \right).
\end{equation}

Conditions \eqref{Bellman_initial},\eqref{Bellman_equation} together with the continuity of $B$ are sufficient to certify optimality of the value $B(t_0,x_0)$ for trajectories starting at $(t_0,x_0)$. Indeed, let $x(t)$, $t \in [t_0,t_f]$, $(t_f,x(t_f)) \in M_f$, $x(t_0) = x_0$ be an admissible trajectory with control $u(t)$. Then $u(t) \in U_{x(t)}$ for almost all $t \in [t_0,t_f]$, because $x(t) \in X$, and the value of the objective functional on this trajectory is given by
\begin{align*}
{\cal I}(x(\cdot)) &= \int_{t_0}^{t_f}L(x(t),u(t))\,dt \leq -\int_{t_0}^{t_f} \frac{\partial B}{\partial t} + \frac{\partial B}{\partial x}g(x(t),u(t)) \, dt \\ =& -\int_{t_0}^{t_f} \frac{dB(t,x(t))}{dt}\,dt = -(B(t_f,x(t_f)) - B(t_0,x(t_0))) = B(t_0,x_0).
\end{align*}
Here the inequality holds by virtue of \eqref{Bellman_equation}, and the last equality by virtue of \eqref{Bellman_initial}. Thus the value of the objective on any admissible trajectory cannot exceed $B(t_0,x_0)$. On the other hand, this value can be achieved if a control $u$ is chosen a.e.\ which realizes the maximum in \eqref{Bellman_equation}, because in this case the inequality turns into an equality.

The Bellman function provides the maximal value of the objective functional for a fixed initial point $(t_0,x_0)$. In order to find the maximal value under the constraint $(t_i,x(t_i)) \in M_i$ one has to maximize $B(t,x)$ over $M_i$. If the objective has to be minimized, then the maximum in \eqref{Bellman_equation} has to be a minimum and we have to minimize $B$ over $M_i$.

In the following sections we prove the theorems by converting them into optimization problems as above, then providing explicit expressions for the corresponding Bellman functions and maximizing (minimizing) them over the manifold of initial points. The proof hence consists in showing continuity, relations \eqref{Bellman_initial} and \eqref{Bellman_equation}, and maximizing (minimizing) over the set of initial points. We shall also sketch how one can arrive at the solution of the problem by solving the Euler-Lagrange equation or applying the PMP without going into much detail.

\section{Proof of Theorem \ref{thm:C_libre}} \label{sec:gen_ineq}

As in \cite{Tholozan17}, is it sufficient to consider the case $n = 1$. Let $\Omega \subset \mathbb RP^1$ be a proper open line segment, i.e., such that its complement is also a line segment. Without loss of generality we assume that the cone $K_{\Omega} \subset \mathbb R^2$ over the closure of $\Omega$ is the positive orthant $\mathbb R_+^2$. Following \cite{Tholozan17}, we parameterize $\Omega$ by a variable $t \in \mathbb R$, in a way such that the point $(e^t,e^{-t}) \in \mathbb R_+^2$ projects to the corresponding point in $\Omega$. Then the Hilbert distance between $t_i,t_f$ is given by $d^H(t_i,t_f) = |t_i-t_f|$.

A centro-affine lift of $\Omega$ into the interior of $\mathbb R_+^2$ is given by a curve $f: t \mapsto e^{\alpha(t)}\cdot(e^t,e^{-t})$. We suppose the function $\alpha(t)$ to be of class $C^2(\mathbb R)$. By virtue of \cite[Lemma 2.1]{Tholozan17} the centro-affine fundamental form of the immersion $f$ is given by $h = \ddot\alpha - \dot\alpha^2 + 1$. Hence locally strong convexity of $f$ is equivalent to the differential inequality $\ddot\alpha > \dot\alpha^2 - 1$. The first derivative is bounded by $|\dot\alpha| < 1$ (\cite[Prop.~2.2]{Tholozan17} with $C = +\infty$). 

We now maximize the Riemannian distance between two given points $t_i < t_f$. This distance is given by the integral
\begin{equation} \label{Riemann_distance}
l^R(t_i,t_f) = \int_{t_i}^{t_f} \sqrt{\ddot\alpha - \dot\alpha^2 + 1}\,dt.
\end{equation}
We shall prove the following estimate.

{\lemma \label{lem:gen_C} Set $d = t_f - t_i$. Then $l^R(t_i,t_f) < \log\left( e^d + \sqrt{e^{2d} - 1} \right)$, and this estimate is sharp. }

The proof is by presenting an explicit Bellman function for the variational problem under consideration. Before proceeding to the proof, we shall give some clues how to arrive at this expression.

Let us cast the problem as an optimal control problem as in Section \ref{sec:Bellman}. Set $x = \dot\alpha$. The extremals for the functional $\int L(x,\dot x)\,dt$ are given by the solutions of the Euler-Lagrange equation $\frac{d}{dt}\frac{\partial L}{\partial\dot x} = \frac{\partial L}{\partial x}$ \cite{Gelfand63}. Inserting $L(x,\dot x) = \sqrt{\dot x - x^2 + 1}$, we obtain the second order ODE $\ddot x = 2x(3\dot x - 2x^2 + 2)$. This ODE has in particular the solutions
\begin{equation} \label{sol_x}
x = -\frac{e^{4t}c^2 + 2(e^{4t} - 2)c + e^{4t}}{(e^{2t} + (e^{2t} - 2)c)(e^{2t}c + e^{2t} - 2)},
\end{equation}
where $c \in (-1,1)$ is an integration constant. The Bellman function is then constructed by computing the value of the cost function on these trajectories. 

\begin{proof}
By invariance with respect to translations of the variable $t$ we may set $t_f = 0$, $t_i = -d$. The set of final points $(t,x)$ is hence given by $M_f = \{ 0 \} \times [-1,1]$, the set of initial points by $M_i = \{ -d \} \times [-1,1]$, the admissible set by $X = [-1,1]$. The dynamics is given by $\dot x = g(x,u) = u \geq x^2-1$, and the cost to maximize by the integral of $L(x,u) = \sqrt{u-x^2+1}$. Consider the function
\[ B(t,x) = \frac12\log\left(\sqrt{(1 - e^{2t})(1 - x)(1 + e^{2t} - x(1 - e^{2t}))} - x(1 - e^{2t}) + 1\right) - t
\]
on $\mathbb R_- \times [-1,1]$. We have $B(0,x) \equiv 0$, and $B$ satisfies \eqref{Bellman_initial}. By the inequality between arithmetic and geometric mean we have
\begin{align*}
\frac{\partial B}{\partial x}u + \frac{\partial B}{\partial t} &= -\frac{(1 - x)(1 + e^{2t} - x(1 - e^{2t})) + (u - x^2 + 1)(1 - e^{2t})}{2\sqrt{(1 - e^{2t})(1 - x)(1 + e^{2t} - x(1 - e^{2t}))}} \\ \leq& -\frac{\sqrt{(1 - x)(1 + e^{2t} - x(1 - e^{2t}))(u - x^2 + 1)(1 - e^{2t})}}{\sqrt{(1 - e^{2t})(1 - x)(1 + e^{2t} - x(1 - e^{2t}))}} = -\sqrt{u - x^2 + 1} = -L(x,u)
\end{align*}
for every $u \geq x^2 - 1$. Here equality is achieved at $(1 - x)(1 + e^{2t} - x(1 - e^{2t})) = (u - x^2 + 1)(1 - e^{2t})$, or equivalently
\[ u = \frac{(1 - x)(1 + e^{2t} - x(1 - e^{2t}))}{1 - e^{2t}} - (1 - x^2) = \frac{2(1 - x)(e^{2t} - x(1 - e^{2t}))}{1 - e^{2t}}.
\]
This proves \eqref{Bellman_equation}, and $B$ is indeed the Bellman function.

We now have to maximize the Bellman function over the initial set $M_i$. In view of
\[ \frac{\partial B}{\partial x} = -\frac{\frac{e^{2t}(1-e^{2t})}{\sqrt{(1 - e^{2t})(1 - x)(1 + e^{2t} - x(1 - e^{2t}))}} + 1 - e^{2t}}{2\left(\sqrt{(1 - e^{2t})(1 - x)(1 + e^{2t} - x(1 - e^{2t}))} - x(1 - e^{2t}) + 1\right)} < 0
\]
for $t < 0$ we obtain that the maximum is achieved at $(t,x) = (-d,-1)$ and given by $\log(1 + \sqrt{1 - e^{-2d}}) + d = \log(e^d + \sqrt{e^{2d} - 1})$.

The trajectory maximizing the objective function is hence obtained by integrating with initial value $x(-d) = -1$ and with $\dot x$ equal to the value of $u$ given above. It is easily verified that this yields the solution \eqref{sol_x} with $e^{-2d} = \frac{4c}{(1+c)^2}$ and the solution satisfies $x(0) = 1$.

We have proven that no centro-affine immersion $f$ can yield a length $l^R(t_i,t_f)$ exceeding the bound in the lemma. Moreover, the bound cannot be attained, because on the optimal trajectory yielding this value we have $x(t_i) = -1$, $x(t_f) = 1$, but a valid immersion satisfies $|x| < 1$. Extend the optimal trajectory by setting $x = -1$ for $t < t_i$ and $x = 1$ for $t > t_f$. Although the resulting function $\alpha(t)$ is not $C^2$, because its second derivatives are discontinuous at $t = t_i$ and $t = t_f$, it can be approximated by $C^2$ functions which satisfy $|\dot\alpha| < 1$ and $\ddot\alpha > \dot\alpha^2 - 1$ everywhere and yield objective values arbitrarily close to the bound in the lemma.
\end{proof}

Let us now return to Theorem \ref{thm:C_libre}. Let $l$ be the projective line passing through the points $a,b \in \Omega$, and $L \subset \mathbb R^{n+1}$ the two-dimensional linear subspace over $l$. Then the centro-affine metric $h_l$ of the immersion $f|_{\Omega \cap l}$ into $L$ is given by the restriction of the centro-affine metric $h$ to $\Omega \cap l$, because the position vector field on $\Omega \cap l$ is contained in $L$. Therefore the Riemannian length $l^R(a,b)$ of the line segment between $a$ and $b$ is equal in both metrics $h$ and $h_l$. Moreover, by definition the Hilbert distance $d^H(a,b)$ is equal in $\Omega$ and in $\Omega \cap l$. The assertion of Lemma \ref{lem:gen_C} is an inequality between these quantities as defined by the immersion of $\Omega \cap l$ into $L \sim \mathbb R^2$. But then the lemma proves also the second inequality in Theorem \ref{thm:C_libre}, which is between the same quantities defined by the immersion $f$ of $\Omega$ into $\mathbb R^{n+1}$, and shows that it cannot be improved.

The first inequality follows from the fact that the geodesic distance between two points in a Riemannian manifold is never exceeding the length of any curve between these two points. For $n = 1$ it turns into an equality, because the straight line segment is the only path linking the two points.

The third inequality in Theorem \ref{thm:C_libre} follows from the relation $\sqrt{e^{2d}-1} < e^d$. Moreover, it is sharp, because
\[ \lim_{d \to \infty} \left( e^d + \sqrt{e^{2d}-1} - d \right) = \log 2.
\]

This proves Theorem \ref{thm:C_libre}.

\section{Proof of Theorem \ref{thm:C_bound}} \label{sec:bound_ineq}

Let us again set $n = 1$ and assume the notations at the beginning of Section \ref{sec:gen_ineq}. However, this time we consider immersions $f: \Omega \to \mathbb R_+^2$ defined by functions $\alpha(t)$ of class $C^3$.

Let us compute the cubic form $C$ of the immersion $f$. It is given by the derivative $\nabla h$, where $h$ is the affine metric and $\nabla$ is the induced affine connection. The latter is defined by the decomposition of the canonical flat affine connection $D$ of the ambient space $\mathbb R^2$ into a tangential part $\nabla$ and a transversal part $h \cdot f$ \cite[p.~28]{NomizuSasaki}. Let $X = \dot f = \dot\alpha f + Jf$, where $J = \diag(1,-1)$, be the basis tangent vector field. Then we have the decomposition
\[ D_XX = \ddot f = \ddot\alpha f + \dot\alpha\dot f + J\dot f = \ddot\alpha f + \dot\alpha^2 f + 2\dot\alpha Jf + f = (\ddot\alpha - \dot\alpha^2 + 1)f + 2\dot\alpha(\dot\alpha f + Jf) = h(X,X) \cdot f + 2\dot\alpha \cdot X.
\]
Therefore $\nabla_XX = 2\dot\alpha \cdot X$. Hence we obtain
\[ C(X,X,X) = Xh(X,X) - 2h(\nabla_XX,X) = \frac{d}{dt}(\ddot\alpha - \dot\alpha^2 + 1) - 4\dot\alpha(\ddot\alpha - \dot\alpha^2 + 1) = \dddot\alpha - 6\dot\alpha\ddot\alpha + 4\dot\alpha^3 - 4\dot\alpha.
\]

Condition \eqref{Cbound} can then be written as
\begin{equation} \label{Cbound2}
|\dddot\alpha - 6\dot\alpha\ddot\alpha + 4\dot\alpha^3 - 4\dot\alpha| \leq 2\gamma(\ddot\alpha - \dot\alpha^2 + 1)^{3/2}.
\end{equation}
We shall now show that functions $\alpha: \mathbb R \to \mathbb R$ satisfying this differential inequality must satisfy certain bounds on the values of their derivatives $\ddot\alpha,\dot\alpha$.

{\lemma \label{lem:h_upper_bound} Let $\alpha: \mathbb R \to \mathbb R$ be a $C^2$ function satisfying \eqref{Cbound2} a.e.\ for some $\gamma \geq 0$. Then for every $t \in \mathbb R$ we have
\[ \sqrt{\ddot\alpha - \dot\alpha^2 + 1} \leq \mu(1 - |\dot\alpha|),
\]
where $\mu \geq 1$ depends on $\gamma$ by the relations $\mu = \frac{\gamma}{2} + \sqrt{1 + \frac{\gamma^2}{4}}$, $\gamma = \frac{\mu^2-1}{\mu}$. }

\begin{proof}
Set $h = \ddot\alpha - \dot\alpha^2 + 1$. Since $h \geq 0$, we have $\ddot\alpha \geq \dot\alpha^2 - 1$. Therefore, if for some $t_0$ we have $\dot\alpha(t_0) > 1$, then $\dot\alpha(t)$ blows up before $t$ reaches $+\infty$. Likewise, if $\dot\alpha(t_0) < -1$, then $\dot\alpha(t)$ blows up before $t$ reaches $-\infty$. Thus $|\dot\alpha| \leq 1$.

Relation \eqref{Cbound2} can be written as
\begin{equation} \label{Cbound3}
-2\gamma h^{3/2} + 4\dot\alpha h \leq \dot h \leq 2\gamma h^{3/2} + 4\dot\alpha h.
\end{equation}
Thus either $h \equiv 0$, in which case the assertion of the lemma is obvious, or $h > 0$ everywhere. In the latter case, define functions $\xi_{\pm} = \dot\alpha \pm \mu^{-1}\sqrt{h}$. By \eqref{Cbound3} we have
\[ \dot\xi_{\pm} = \ddot\alpha \pm \mu^{-1}\frac{\dot h}{2\sqrt{h}} \geq \ddot\alpha + \mu^{-1}\frac{-2\gamma h^{3/2} \pm 4\dot\alpha h}{2\sqrt{h}} = h + \dot\alpha^2 - 1 + \mu^{-1}(-\gamma h \pm 2\dot\alpha\sqrt{h}) = \xi_{\pm}^2 - 1.
\]
Here we used that $1 - \mu^{-1}\gamma = \mu^{-2}$.

As above we have $|\xi_{\pm}| \leq 1$ and hence $\dot\alpha + \mu^{-1}\sqrt{h} \leq 1$, $-1 \leq \dot\alpha - \mu^{-1}\sqrt{h}$. It follows that $\sqrt{h} \leq \mu(1 \pm \dot\alpha)$.
\end{proof}

{\lemma \label{lem:h_lower_bound} Suppose in addition to the conditions in Lemma \ref{lem:h_upper_bound} that the function $\alpha$ defines a non-degenerate centro-affine immersion $f: t \mapsto e^{\alpha(t)}\cdot(e^t,e^{-t})$ which is asymptotic to the boundary of the orthant $\mathbb R_+^2$. Then for every $t \in \mathbb R$ we have that
\[ \mu^{-1}(1+|\dot\alpha|) \leq \sqrt{\ddot\alpha - \dot\alpha^2 + 1}.
\] }

\begin{proof}
Let $h = \ddot\alpha - \dot\alpha^2 + 1 > 0$ be the affine metric. Define functions $\psi_{\pm} = \dot\alpha \pm \mu\sqrt{h}$. By \eqref{Cbound3} we then have
\[ \dot\psi_{\pm} = \ddot\alpha \pm \mu\frac{\dot h}{2\sqrt{h}} \leq \ddot\alpha + \mu\frac{2\gamma h^{3/2} \pm 4\dot\alpha h}{2\sqrt{h}} = h + \dot\alpha^2 - 1 + \mu(\gamma h \pm 2\dot\alpha\sqrt{h}) = \psi_{\pm}^2 - 1.
\]
Here we used that $1 + \mu\gamma = \mu^2$.

Now suppose for the sake of contradiction that for some $t_0$ we have $\psi_+(t_0) < 1$. Choose $\eta_0$ such that $\psi_+(t_0) < -\tanh(t_0 + \eta_0)$. Then we have $\kappa(t) = \psi_+(t) + \tanh(t + \eta_0) < 0$ for all $t > t_0$. Indeed, let $t_1 > t_0$ be the smallest point such that $\kappa(t_1) = 0$. For all $t < t_1$ close enough to $t_1$ we have $\kappa(t) < 0$, $\psi_+(t) \in (-1,1)$ and
\[ \frac{d\kappa}{dt} \leq \psi_+(t)^2 - 1 + \frac{1}{\cosh(t + \eta_0)^2} = (\psi_+(t) - \tanh(t + \eta_0))\cdot\kappa(t) < -2\kappa(t).
\]
But then $\kappa$ cannot reach zero in finite time, a contradiction. Hence such a point $t_1$ cannot exist. We obtain $\dot\alpha(t) + \mu\sqrt{h(t)} < -\tanh(t + \eta_0)$ for all $t \geq t_0$ and by integration
\[ 0 < \mu\int_{t_0}^{t}\sqrt{h(s)}\,ds < \left[ -\alpha(s) - \log\cosh(s + \eta_0) \right]_{t_0}^t = \left[ -(s+\alpha(s)) - \log\frac{e^{\eta_0} + e^{-(2s + \eta_0)}}{2} \right]_{t_0}^t.
\]
But $t+\alpha(t) \to +\infty$ for $t \to +\infty$, because $e^{t+\alpha}$ is the first coordinate of the immersion $f$, which by assumption is asymptotic to the boundary of $\mathbb R_+^2$. Therefore the right-most expression tends to $-\infty$ as $t \to +\infty$, a contradiction. As a consequence, we have $\psi_+ \geq 1$ and $\mu\sqrt{h} \geq 1 - \dot\alpha$.

The inequality $\mu\sqrt{h} \geq 1 + \dot\alpha$ is proven similarly by deducing a contradiction from the assumption $\psi_-(t_0) > -1$ for some $t_0$.
\end{proof}

{\corollary \label{cor:infinitesimal_length} Assume the conditions of Lemmas \ref{lem:h_upper_bound}, \ref{lem:h_lower_bound}, and let $t_i < t_f$ be points in $\Omega = \mathbb R$. Then the Riemannian length $l^R(t_i,t_f)$ is bounded by $\mu^{-1}d^H(t_i,t_f) \leq l^R(t_i,t_f) \leq \mu d^H(t_i,t_f)$. }

\begin{proof}
By Lemmas \ref{lem:h_upper_bound}, \ref{lem:h_lower_bound} we have $\sqrt{\ddot\alpha - \dot\alpha^2 + 1} \in [\mu^{-1},\mu]$. The claim now follows from \eqref{Riemann_distance}.
\end{proof}

{\corollary Assume the conditions of Lemmas \ref{lem:h_upper_bound}, \ref{lem:h_lower_bound}. If $\gamma = 0$, or equivalently $\mu = 1$, then $\alpha \equiv const$, and the image of $f$ is a hyperbola. In this case the Riemannian length $l^R(t_i,t_f)$ coincides with the Hilbert distance $d^H(t_i,t_f)$. }

\begin{proof}
From $\mu = 1$ we have $1 + |\dot\alpha| \leq 1 - |\dot\alpha|$ by virtue of Lemmas \ref{lem:h_upper_bound}, \ref{lem:h_lower_bound}. We get $\alpha \equiv const$, $h \equiv 1$, and $l^R(t_i,t_f) = |t_i - t_f| = d^H(t_i,t_f)$.
\end{proof}

We may write the problem of maximizing the Riemannian length $l^R(t_i,t_f)$ under constraint \eqref{Cbound3} as an optimal control problem. Introduce variables $x = \dot\alpha$, $y = \sqrt{h}$. Then the dynamics of the system can be written as
\begin{equation} \label{control_dynamics}
\dot x = y^2 + x^2 - 1, \quad \dot y = 2xy + u\gamma y^2, \qquad u \in [-1,1].
\end{equation}
Here the first equation comes from the definition of $y$, and the second equation is equivalent to \eqref{Cbound3}. The variable $u$ is a scalar control. The objective is to maximize
\begin{equation} \label{control_cost}
l^R(t_i,t_f) = \int_{t_i}^{t_f} y(t)\,dt \to \sup.
\end{equation}
In addition we have the state constraints
\begin{equation} \label{state_constraints}
\mu^{-1}(1+|x|) \leq y \leq \mu(1-|x|)
\end{equation}
from Lemmas \ref{lem:h_upper_bound}, \ref{lem:h_lower_bound}. Relations \eqref{state_constraints} define the feasible set $X \subset \mathbb R^2$.

It is easily checked that for a constant control $u$ the trajectories of system \eqref{control_dynamics} are given by the level curves of the first integral
\[ I_u = \frac{(\mu_u^2 + 1)y}{\mu_u(y^2 - x^2 + 1) - (\mu_u^2 - 1)xy},
\]
where $\mu_u = \frac{u\gamma}{2} + \sqrt{1 + \frac{u^2\gamma^2}{4}}$.

\medskip

Theorem \ref{thm:C_bound} will be proven by means of Lemmas \ref{lem:defB}, \ref{lem:maximalB} below. The proofs of the lemmas give no clue how the Bellman functions in these lemmas have been obtained, however. Before we state the lemmas, we shall therefore sketch how one can arrive at these functions by optimal control techniques.

\begin{figure}
\centering
\includegraphics[width=16.18cm,height=4.18cm]{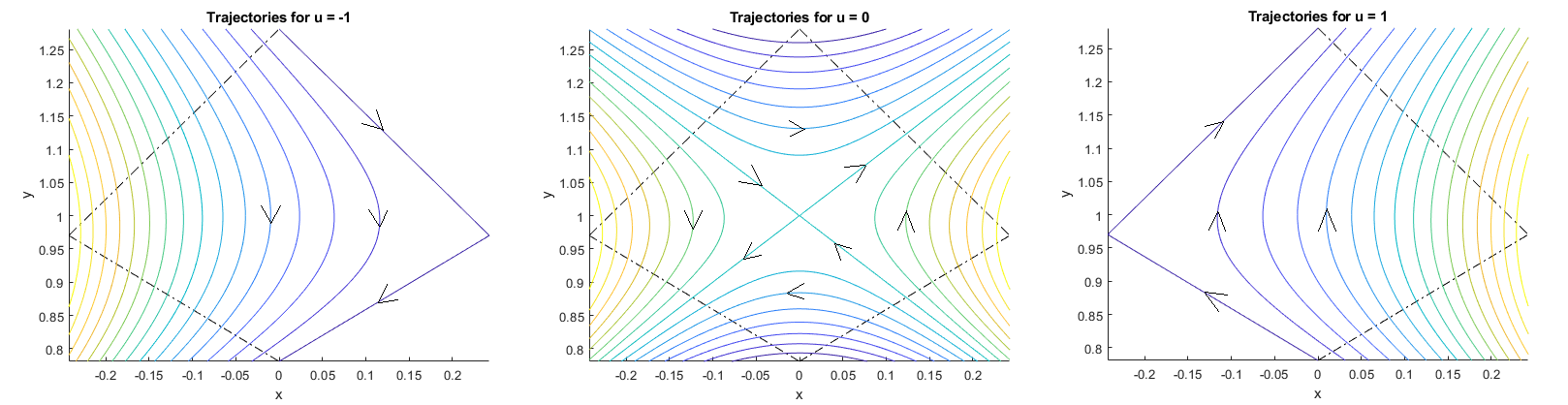}
\caption{Trajectories of system \eqref{control_dynamics} for different constant values of $u$. The parameter $\gamma$ equals $0.5$. The dash-dotted lines delimit the feasible region $X$ given by state constraints \eqref{state_constraints}.}
\label{fig:trajectories}
\end{figure}

The boundary segments of $X$ are trajectories of the system with constant control $u = \pm1$ (see Fig.~\ref{fig:trajectories}). Moreover, no feasible trajectory of the system can leave the upper right and the lower left boundary segment after hitting it. Likewise, no feasible trajectory can leave the upper left and the lower right boundary segment after hitting it in backward time. If a trajectory of \eqref{control_dynamics} leaves $X$, then it cannot return to $X$ anymore. Therefore, if the initial point $(x(t_i),y(t_i))$ and the terminal point $(x(t_f),y(t_f))$ of a trajectory of \eqref{control_dynamics} satisfy \eqref{state_constraints}, then all intermediate points do so too. Thus we do not need to take the state constraints into account when considering the first order optimality conditions for a trajectory with fixed initial and terminal point. 

Let us apply the PMP to control problem \eqref{control_dynamics},\eqref{control_cost} with fixed initial and terminal point \cite{PBGM62}. According to this principle, if $(x(t),y(t))$, $t \in [t_i,t_f]$, is a maximizer of \eqref{control_cost} under the additional constraints $(x(t_i),y(t_i)) = (x_i,y_i)$, $(x(t_f),y(t_f)) = (x_f,y_f)$, then there exist a nonnegative constant $\lambda$ and differentiable functions $\psi(t),\phi(t)$ (the so-called adjoint variables), not all equal to zero, such that at every $t \in [t_i,t_f]$ the control $u$ maximizes the \emph{Pontryagin function}
\[ {\cal H} = 2\lambda y + \psi\cdot(y^2 + x^2 - 1) + \phi\cdot(2xy + u\gamma y^2),
\]
and the adjoint variables are solutions of the differential equations
\begin{equation} \label{adjoint_dynamics}
\dot\psi = -\frac{\partial{\cal H}}{\partial x} = -2(x\psi + y\phi), \quad \dot\phi = -\frac{\partial{\cal H}}{\partial y} = -2(\lambda + y\psi + x\phi + u\gamma y\phi).
\end{equation}
It follows that the control $u$ is given by the sign of $\phi$ whenever this variable does not vanish.

If the partial derivative $\frac{\partial{\cal H}}{\partial u}$ vanishes identically on some interval, then the corresponding trajectory is called a \emph{singular arc}. Let us compute the corresponding control $u$. If $\phi \equiv 0$, then also $\dot\phi \equiv 0$, which entails $\lambda + y\psi \equiv 0$. If $\lambda = 0$, then also $\psi = 0$, and by \eqref{adjoint_dynamics} the adjoint variables vanish identically. Thus the presence of a singular arc entails $\lambda > 0$. Differentiating further, we obtain $\dot y\psi + y\dot\psi = u\gamma y^2\psi = -\lambda u\gamma y \equiv 0$, entailing $u \equiv 0$.

On optimal trajectories the control $u$ is hence piece-wise constant and taking values in $\{-1,0,+1\}$. Any optimal trajectory can therefore be assembled from the trajectories depicted in Fig.~\ref{fig:trajectories}. If some trajectory is optimal in the larger class of trajectories with free end-points, then the necessary optimality conditions obtained for trajectories with fixed end-points still apply. After some calculations one obtains the following solution.

{\lemma \label{lem:defB} Let $t_i < t_f$, $T = t_f - t_i$, $\gamma > 0$, $\mu = \frac{\gamma}{2}+\sqrt{1+\frac{\gamma^2}{4}}$,
\begin{equation} \label{abcd}
a_{\pm} = \mu y \pm (1-x),\ b_{\pm} = \mu y \pm (1+x),\ c_{\pm} = \mu(1-x) \pm y,\ d_{\pm} = \mu(1+x) \pm y,
\end{equation}
and consider control problem \eqref{control_dynamics}--\eqref{state_constraints} with fixed initial point $(x(t_i),y(t_i)) = (x_i,y_i)$ and free terminal point $(x(t_f),y(t_f)) \in X$. Then the optimal value of the problem is given by $B(-T,x_i,y_i)$, where the function $B: \mathbb R_- \times X \to \mathbb R_+$ is defined as follows.

Let $(x,y) \in X$ and $t \in \mathbb R_-$. If
\[ -t \leq t_{+1}(x,y) = \left\{ \begin{array}{rcl} \frac12\log\left( 1 + \frac{2c_-}{a_-c_+} \right),&\quad& a_- > 0, \\ +\infty, && a_- = 0, \end{array} \right.
\]
then
\[ B(t,x,y) = \frac{\mu}{\mu^2+1}\log\frac{c_+e^{-2t} + d_-}{\mu(-a_-e^{-2t} + b_+)}.
\]
If
\[ t_{+1}(x,y) < -t \leq \hat t(x,y) = \left\{ \begin{array}{rcl} \frac12\log\frac{b_+}{a_-},&\quad& a_- > 0, \\ +\infty, && a_- = 0, \end{array} \right.
\]
then
\begin{align*}
B(t,x,y) &= \frac{\mu}{\mu^2+1}\log\left( \frac{e^{-2t}(\mu^2 + 1)c_+^2}{8\mu y} + \right. \\ & + \left. \frac{c_+(\mu^4xy - \mu^4y + \mu^3x^2 - \mu^3y^2 - \mu^3 + 6\mu^2y + \mu x^2 - \mu y^2 - \mu - xy - y)}{8\mu ya_-}\right).
\end{align*}
If
\[ \hat t(x,y) < -t \leq t^*(x,y) = \left\{ \begin{array}{rcl} \frac12\log\frac{1 - x^2 + y^2}{(y + 1 - x)(x + y - 1)},&\quad& x + y > 1, \\ +\infty, && x + y \leq 1, \end{array} \right.
\]
then
\begin{align*}
\lefteqn{B(t,x,y) = \frac12\log\frac{W + a_-c_+e^{-2t} + \mu^2xy + \mu x^2 - \mu y^2 - \mu - xy}{(\mu+1)^2y} +} \\ & + \frac{\mu}{\mu^2+1}\log\frac{2\mu^2\left[W(1+\mu^2) - 2\mu(a_-c_+e^{-2t} + \mu^2xy + \mu x^2 - \mu y^2 - \mu - xy)\right]}{(\mu^2-1)^2a_-(a_-e^{-2t} - b_+)}
\end{align*}
with
\[ W = \sqrt{(c_+a_-e^{-2t} + \mu^2xy + \mu x^2 - \mu y^2 - \mu - xy)^2 - (\mu^2y - y)^2}.
\]
If $t^*(x,y) < -t$, then
\begin{align*}
\lefteqn{B(t,x,y) = \frac12\log\frac{W + c_-a_+e^{-2t} + \mu^2xy - \mu x^2 + \mu y^2 + \mu - xy}{(\mu+1)^2y} +} \\ & + \frac{\mu}{\mu^2+1}\log\frac{2\mu^2\left[(\mu^2 + 1)W - 2\mu(c_-a_+e^{-2t} + \mu^2xy - \mu x^2 + \mu y^2 + \mu - xy)\right]}{(\mu^2-1)^2(c_-e^{-2t} + d_+)c_-}
\end{align*}
with
\[ W = \sqrt{(c_-a_+e^{-2t} + \mu^2xy - \mu x^2 + \mu y^2 + \mu - xy)^2 - (\mu^2y - y)^2}.
\]}

\begin{proof}
As in the proof of Lemma \ref{lem:gen_C} we may set $t_f = 0$, and hence $M_f = \{0\} \times X$. The dynamics of the system is given by $g(x,y,u) = (g_1,g_2) = (y^2+x^2-1,2xy+u\gamma y^2)$, the objective by the integral of $L(x,y,u) = y$.

The proof mainly consists of showing a series of inequalities. Transforming the corresponding expressions involves many calculations, which cannot all be included, but are rather straightforward. We shall concentrate on certifying the inequalities by presenting the corresponding expressions in a form suitable for immediately recognizing their nonnegativity or non-positivity.

Let us show that $B(t,x,y)$ is the Bellman function of the problem. This involves showing that $B$ is well-defined, continuous, satisfies \eqref{Bellman_initial},\eqref{Bellman_equation}, and the optimal control $\hat u$ yielding the maximum in \eqref{Bellman_equation} defines a feasible trajectory. Denote the four expressions defining $B$ in the lemma by $B_I,B_{II},B_{III},B_{IV}$. Denote the expressions for $W$ in the lemma by $W_{III},W_{IV}$. Recall that $y > 0$, $1 \pm x > 0$, $\mu > 1$. By \eqref{state_constraints} the quantities $a_{\pm},b_{\pm},c_{\pm},d_{\pm}$ are positive in the interior of $X$.

\medskip

\emph{Consistency:} For fixed $(x,y) \in X$ the values $0,t_{+1}(x,y),\hat t(x,y),t^*(x,y)$ form an increasing sequence. Indeed, clearly $t_{+1}(x,y) \geq 0$. Straightforward calculation yields
\[ \hat t(x,y) - t_{+1}(x,y) = -\frac12\log\left( 1 - \frac{4y}{b_+c_+} \right) > 0,
\]
and for $x+y > 1$ we get
\[ t^*(x,y) - \hat t(x,y) = \frac12\log\left( 1 + \frac{2yc_-}{b_+(y + 1 - x)(x + y - 1)} \right) \geq 0.
\]

According to the decomposition $\alpha^2 - \beta^2 = (\alpha - \beta)(\alpha + \beta)$ the expression $W_{III}$ is the square root of the product of two functions which are both linear in $e^{-2t}$ with nonnegative leading coefficient $c_+a_-$. Inserting $t = -\hat t$ into these linear factors, one easily calculates that these evaluate to $2y$ and $2\mu^2 y$, respectively. For general $-t \geq \hat t$ these factors are hence nonnegative, and $W_{III}$ is well-defined. Likewise, the linear factors involved in the definition of $W_{IV}$ have leading coefficient $a_+c_- \geq 0$. Inserting $t = -t^*$ yields the positive values $\frac{2y(1 - x)^2(\mu^2 - 1)}{(y + 1 - x)(x + y - 1)}$, $\frac{2y^3(\mu^2 - 1)}{(y + 1 - x)(x + y - 1)}$. Therefore $W_{IV}$ is well-defined for $-t \geq t^*$.

\medskip

\emph{Continuity:} Inserting $t = -t_{+1}(x,y)$ into $B_I$ and $B_{II}$, we obtain the same expression
\[ \frac{\mu}{\mu^2+1} \log\frac{(\mu^2-1)c_+}{2\mu a_-}.
\]
Inserting $t = -\hat t(x,y)$ into $B_{II}$ and $B_{III}$, we obtain the same expression
\[ \frac{\mu}{\mu^2+1}\log\frac{\mu c_+}{a_-}.
\]
Inserting $t = -t^*(x,y)$ into $B_{III}$ and $B_{IV}$, we obtain for $x+y > 1$ the same expression
\[ \frac12\log\frac{(\mu-1)(y - x + 1)}{(\mu+1)(y + x - 1)} + \frac{\mu}{\mu^2+1}\log\frac{2\mu^2}{\mu^2 - 1}.
\]
Hence $B$ is continuous.

\medskip

\emph{Initial value:} Inserting $t = 0$ into $B_I$, we obtain $B(0,x,y) \equiv 0$, which proves \eqref{Bellman_initial}.

\medskip

\emph{Bellman inequality:} Let us show \eqref{Bellman_equation}, i.e., that for all $u \in [-1,1]$ we have
\begin{equation} \label{Bellman_ineq}
\frac{d}{dt}B(t,x(t),y(t)) + y = \frac{\partial B}{\partial t} + \frac{\partial B}{\partial x} \cdot g_1 + \frac{\partial B}{\partial y} \cdot g_2 + y \leq 0,
\end{equation}
with equality if
\[ u = \hat u(t,x,y) = \left\{ \begin{array}{rcl} +1,&\quad& -t < t^*(x,y)\ \mbox{and}\ c_- > 0, \\ 0,&& -t = t^*(x,y)\ \mbox{and}\ c_- > 0,\\ -1,&& -t > t^*(x,y)\ \mbox{or}\ c_- = 0. \end{array} \right.
\]

For $0 \leq -t < t_{+1}$ we have
\[ \frac{dB_I}{dt} + y = -\frac{y^2(\mu^2 - 1)(e^{-2t} - 1)(1 - u)(e^{-2t}(1-x) + 1 + x)c_+}{(a_-c_+(e^{2t_{+1}}-e^{-2t}) + 4y)(c_+e^{-2t} + d_-)} \leq 0
\]
with equality if $u = 1$. Note that if $c_- = 0$, then $t_{+1} = 0$, and hence the inequality $-t < t_{+1}$ cannot hold.

For $t_{+1} \leq -t \leq \hat t$ we have
\[ \frac{dB_{II}}{dt} + y = -\frac{y(\mu^2 - 1)(1 - u)c_-(c_+a_-^2(e^{2\hat t}-e^{-2t}) + 2\mu yc_-)}{a_-c_+\left[(\mu^2 + 1)a_-c_+(e^{-2t}-e^{2t_{+1}}) + 4y(\mu^2 - 1)\right]} \leq 0,
\]
with equality if $u = 1$ or $c_- = 0$.

For $\hat t < -t \leq t^*$ we have
\[ \frac{dB_{III}}{dt} + y = -\frac{y(\mu^2-1)(1-u)A}{2a_-(e^{-2t}-e^{2\hat t})c_+(c_+a_-(e^{-2t}-e^{2\hat t}) + 2y(\mu^2 + 1))},
\]
where $A = \left((y^2 + (1 - x)^2)(e^{-2t}-e^{2\hat t}) + \frac{2yc_+}{a_-}\right)W_{III} + A'$ and
\begin{align*}
\lefteqn{\left((y^2 + (1 - x)^2)(e^{-2t}-e^{2\hat t}) + \frac{2yc_+}{a_-}\right)^2W_{III}^2 - (A')^2 = } \\ =& \left((x - y - 1)(x + y - 1)e^{-2t} + y^2 + 1 - x^2\right)^2  c_+a_-(c_+e^{-2t} + d_-)a_-(e^{-2t}-e^{2\hat t}) \geq 0.
\end{align*}
Since the coefficient at $W_{III}$ in $A$ is nonnegative, we obtain $A \geq 0$ and hence $\frac{dB_{III}}{dt} + y \leq 0$, with equality if $u = 1$. Note that if $c_- = 0$, then $\hat t = t^*$, and the condition $\hat t < -t \leq t^*$ cannot hold.

For $-t > t^*$ we have
\[ \frac{dB_{IV}}{dt} + y = -\frac{y(\mu^2-1)(1+u)a_-A}{2c_-(c_-e^{-2t} + d_+)a_+(a_+a_-(e^{-2t}-e^{2\hat t}) + 4\mu y)},
\]
where $A = \left((y^2 + (1 - x)^2)(e^{-2t} - e^{2t^*}) + \frac{4y^2(1 - x)}{(x + y - 1)(y + 1 - x)}\right)W_{IV} + A'$ and
\begin{align*}
\lefteqn{\left((y^2 + (1 - x)^2)(e^{-2t} - e^{2t^*}) + \frac{4y^2(1 - x)}{(x + y - 1)(y + 1 - x)}\right)^2W_{IV}^2 - (A')^2 = } \\ =& (x - y - 1)^2(x + y - 1)^2(e^{-2t} - e^{2t^*})^2c_-a_+(c_-e^{-2t} + d_+) (a_+(e^{-2t}-1) + 2) \geq 0.
\end{align*}
Since the coefficient at $W_{IV}$ in $A$ is nonnegative, we obtain $A \geq 0$ and hence $\frac{dB_{IV}}{dt} + y \leq 0$, with equality if $u = -1$.

\medskip

\emph{Feasibility:} Let us show that application of the optimal control $\hat u$ guarantees that the trajectory does not leave the feasible region $X$, and indeed arrives at the terminal manifold $M_f$. The only boundary segments through which a trajectory can escape $X$ are the upper right and the lower left one. On the upper right segment we have $c_- = 0$ and hence the optimal control is $\hat u = -1$. The trajectory then moves along the boundary segment. On the lower left segment we have  $t^* = +\infty$ and $c_- > 0$. Therefore the optimal control is $\hat u = +1$ and the trajectory again moves along the boundary segment.

\medskip

Thus $B$ is indeed the Bellman function, and the optimal value is achieved by applying the control $u = \hat u$.
%
\end{proof}

\begin{figure}
\centering
\includegraphics[width=16.19cm,height=3.28cm]{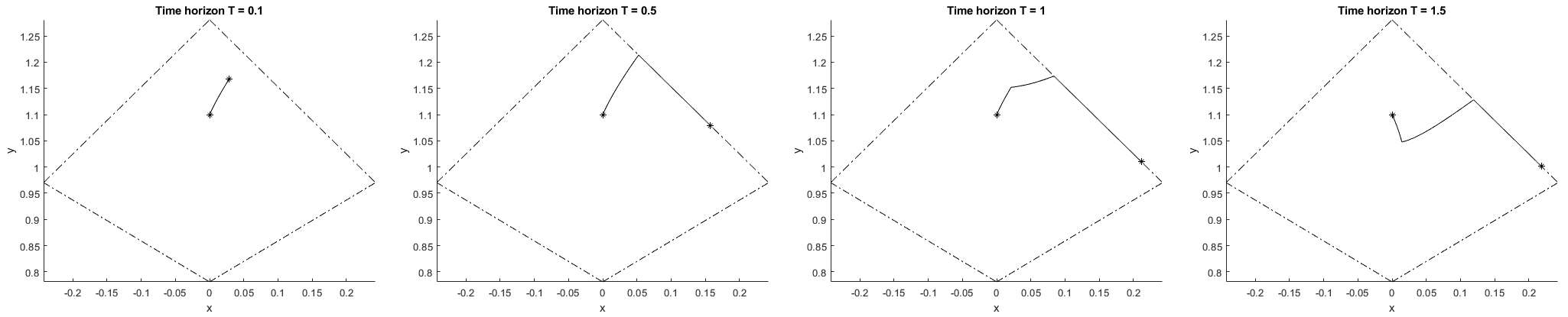}
\caption{Optimal trajectories of system \eqref{control_dynamics} with initial point $(x,y) = (0,1.1)$ for different time horizons $T = t_f - t_i$. The initial and terminal points are marked with stars. The parameter $\gamma$ equals $0.5$. The dash-dotted lines delimit the feasible region $X$.}
\label{fig:optimal_synthesis}
\end{figure}

The corresponding optimal solutions are structured as follows.

If $T \leq t_{+1}(x_i,y_i)$, where $t_{+1}$ is the time needed to reach the upper right boundary segment of the feasible set with control $u \equiv +1$, then the optimal control is given by $\hat u \equiv +1$ on the whole trajectory (see Fig.~\ref{fig:optimal_synthesis}, upper left).

If $t_{+1}(x_i,y_i) < T \leq \hat t(x_i,y_i)$, then on the optimal trajectory the control $u \equiv +1$ is optimal for $t \in [t_i,t_i+t_{+1}(x_i,y_i))$, that is up to the point when the trajectory reaches the boundary of the feasible set. For $t \in (t_i+t_{+1}(x_i,y_i),t_f]$ the control $u \equiv -1$ is optimal and the trajectory moves along the boundary of the feasible set (see Fig.~\ref{fig:optimal_synthesis}, upper right).

If $\hat t(x_i,y_i) < T \leq t^*(x_i,y_i)$, then between the arcs with controls $u = \pm1$ there appears a singular arc with control $u \equiv 0$ (see Fig.~\ref{fig:optimal_synthesis}, lower left).

Finally, for $T > t^*(x_i,y_i)$ the optimal trajectory consists of three arcs, on which the control equals $-1,0,-1$, respectively, with the third arc located on the upper right boundary segment (see Fig.~\ref{fig:optimal_synthesis}, lower right).

\medskip

In order to find the optimal value of problem \eqref{control_dynamics}--\eqref{state_constraints} with free initial and terminal points, we have to maximize the Bellman function $B(-T,x,y)$ over $(x,y) \in X$ for fixed $T$.

{\lemma \label{lem:maximalB} The maximum $\max_{(x,y) \in X}B(-T,x,y)$ is attained at
\[ (x,y) = \left\{ \begin{array}{rcl} \left(-\frac{(\mu^2 - 1)(e^T - 1)}{\mu^2(e^T - 1) + e^T + 1},\frac{2\mu e^T}{\mu^2(e^T - 1) + e^T + 1}\right),&\quad& T \leq \log\frac{\mu^2+1}{\mu^2-1}, \\
\left(- 1 + \frac{e^T}{\mu\sqrt{e^{2T} - 1}},\frac{e^T}{\sqrt{e^{2T} - 1}}\right),&& T \geq \log\frac{\mu^2+1}{\mu^2-1}. \end{array} \right.
\]
The corresponding value of the maximum is given by
\begin{equation} \label{upper_bound_T}
\left\{ \begin{array}{rcl} \frac{2\mu}{\mu^2+1}\log\frac{(\mu^2+1)(e^T - 1) + 2}{2},&\quad& T \leq \log\frac{\mu^2+1}{\mu^2-1}, \\
\log\frac{(\mu-1)(\sqrt{e^T+1} + \sqrt{e^T-1})}{(\mu+1)(\sqrt{e^T+1} - \sqrt{e^T-1})} + \frac{2\mu}{\mu^2+1}\log\frac{2\mu^2}{\mu^2-1},&& T \geq \log\frac{\mu^2+1}{\mu^2-1}. \end{array} \right.
\end{equation}
}

\begin{proof}
We shall parameterize $X$ by the variables $w = \frac{y}{1-x}$, $z = \frac{y}{1+x}$, which both run through the interval $[\mu^{-1},\mu]$. Set $t = -T$. We shall determine the maximum of $B(t,x,y)$ by examining the signs of the derivatives with respect to these variables.

Let us show that $\frac{\partial B}{\partial z} = -\frac{(1-x)(1+x)^2}{2y}\frac{\partial B}{\partial x} + \frac{(1+x)^2}{2}\frac{\partial B}{\partial y} \geq 0$. For $-t \leq t_{+1}(x,y)$ we have
\[ \frac{\partial B_I}{\partial z} = \frac{\mu(x + 1)^2(e^{-2t} - 1)}{a_-(e^{2\hat t}-e^{-2t})(c_+e^{-2t} + d_-)} \geq 0.
\]
For $t_{+1}(x,y) \leq -t \leq \hat t(x,y)$ we have
\[ \frac{\partial B_{II}}{\partial z} = \frac{\mu(e^{-2t} - 1)(x + 1)^2a_-c_+}{2y\left[(\mu^2 + 1)a_-c_+(e^{-2t} - e^{2t_{+1}}) + 4y(\mu^2 - 1)\right]} \geq 0.
\]
For $\hat t(x,y) \leq -t \leq t^*(x,y)$ we have
\[ \frac{\partial B_{III}}{\partial z} = -\frac{(x + 1)^2(e^{-2t} - 1)(W_{III} - 2\mu y)}{4ya_-(e^{2\hat t}-e^{-2t})(c_+e^{-2t} + d_-)} \geq 0.
\]
For $-t \geq t^*(x,y)$ we have
\[ \frac{\partial B_{IV}}{\partial z} =  \frac{(x + 1)^2(e^{-2t} - 1)a_-(W_{IV} + 2\mu y)}{4y(c_-e^{-2t} + d_+)(a_+a_-(e^{-2t}-e^{2\hat t}) + 4\mu y)} \geq 0.
\]
Hence the maximum of $B$ is achieved at $z = \mu$. This corresponds to the upper left boundary segment of $X$.

We now compute the derivative $\frac{\partial B}{\partial w} = \frac{(1-x)^2(1+x)}{2y}\frac{\partial B}{\partial x} + \frac{(1-x)^2}{2}\frac{\partial B}{\partial y}$ on this segment. Note that on this segment $x \in [-\frac{\mu^2-1}{\mu^2+1},0]$, and $e^{2t_{+1}} = \frac{(\mu^2 - 1)(1 + x)}{(\mu^2 + 1)x + \mu^2 - 1}$, $e^{2\hat t} = \frac{(\mu^2 + 1)(1 + x)}{(\mu^2 + 1)x + \mu^2 - 1}$, $e^{2t^*} = \frac{\mu^2(1 + x)^2 + 1 - x^2}{((\mu + 1)x + \mu - 1)(\mu(1 + x) + 1 - x)}$ for $x > -\frac{\mu-1}{\mu+1}$.

For $-t \leq t_{+1}(x,y)$ we have
\[ \frac{\partial B_I}{\partial w} = \frac{\mu e^{-2t}(e^{-2t}-1)(e^{2t_{+1}}-1)(1-x)^2}{2(c_+e^{-2t} + d_-)((1 + x)(e^{-2t}-1) + e^{2t_{+1}}-e^{-2t})} \geq 0.
\]

For $t_{+1}(x,y) \leq -t \leq \hat t(x,y)$ we have
\begin{align*}
\frac{\partial B_{II}}{\partial w} =& -\!\frac{\mu^2(1 - x)^2(e^{2\hat t} - e^{2t_{+1}})\!\left[((\mu^2+1)x + \mu^2 - 1)e^{-t} + (\mu^2 - 1)(1+x)\right]}{2y(1 + x)a_-c_+\left[2\mu^2(e^{-2t} - e^{2t_{+1}}) + (\mu^2 - 1)(e^{2\hat t} - e^{-2t})\right]} \cdot \\
& \cdot \left[((\mu^2 + 1)(e^{-t}-1)+2)x + (\mu^2 - 1)(e^{-t} - 1)\right].
\end{align*}
Hence $\frac{\partial B_{II}}{\partial w} \geq 0$ if $x \leq -\frac{(\mu^2 - 1)(e^{-t} - 1)}{(\mu^2 + 1)(e^{-t}-1)+2}$ and $\frac{\partial B_{II}}{\partial w} \leq 0$ if $x \geq -\frac{(\mu^2 - 1)(e^{-t} - 1)}{(\mu^2 + 1)(e^{-t}-1)+2}$.

For $\hat t(x,y) \leq -t \leq t^*(x,y)$ we have
\begin{align*} 
\frac{\partial B_{III}}{\partial w} =& \frac{2\mu^2(1-x)^2e^{-2t}\beta^2(e^{-2t} - \mu^2(1+x)^2(e^{-2t} - 1))}{ya_-^2c_+(c_+e^{-2t} + d_-)} \cdot \\
&\cdot \frac{1}{\sqrt{(\beta e^{-2t} - x - 1)(\beta e^{-2t} - \mu^2(1+x))} + \mu(1+x)(\beta(e^{-2t}-e^{2\hat t}) + 1)},
\end{align*}
where we wrote $\beta$ for $(\mu^2+1)x + \mu^2 - 1$ for brevity. Hence $\frac{\partial B_{III}}{\partial w} \geq 0$ for $x \leq \frac{e^{-t}}{\mu\sqrt{e^{-2t} - 1}}-1$ and $\frac{\partial B_{III}}{\partial w} \leq 0$ for $x \geq \frac{e^{-t}}{\mu\sqrt{e^{-2t} - 1}}-1$.

For $-t \geq t^*(x,y)$ we have
\[ \frac{\partial B_{IV}}{\partial w} = -\frac{\mu^2(1-x)^2A}{yc_-(c_-e^{-2t} + d_+)a_+(a_+(e^{-2t} - 1) + 2)},
\]
where
\begin{align*}
\lefteqn{A = \left[((\mu^2-1)(x + 1)^2 + 2)(e^{-2t}-1) + 2\right] \cdot} \\ &\cdot \sqrt{\left[((\mu^2+1)x + \mu^2 - 1)e^{-2t} - 1 - x\right]\left[((\mu^2+1)x + \mu^2 - 1)e^{-2t} - \mu^2(1+x)\right]} \\ &+ \mu(x + 1)\left[2x((\mu^2-1)x + \mu^2 + 1)e^{-4t} - (3(\mu^2-1)x^2 + \mu^2(4x + 1) + 1)e^{-2t} + (\mu^2 - 1)(x + 1)^2\right].
\end{align*}
Let us show that $A$ is nonnegative. First we consider the second summand. This is a concave quadratic polynomial in $e^{-2t}$. If we replace $e^{-2t}$ by 0, we obtain the positive value $\mu(\mu^2 - 1)(x + 1)^3$. If we replace $e^{-2t}$ by 1, we obtain the negative value $-2\mu(x+1)$. Therefore for $e^{-2t} \geq 1$ the value of the second term is negative. Since the first summand of $A$ is positive, the difference of the two summands will also be positive. Multiplying $A$ by the difference of the two terms we get rid of the square root, and the resulting expression equals
\[ -x\beta(1 + x - e^{-2t}x)(\beta(e^{-2t} - 1) + 2)\left((\mu^2 - 1)^2(e^{-2t} - 1)^2(x + 1)^4 + 4(x + 1)^2e^{-2t} - 4xe^{-4t}(2+x)\right),
\]
which consists of nonnegative factors. Here we denoted $(\mu^2-1)x + \mu^2 + 1$ by $\beta$. Thus $A \geq 0$ and hence $\frac{\partial B_{IV}}{\partial w} \leq 0$.

It follows that for every fixed $t < 0$ the function $B$ is unimodal on the upper left boundary segment of $X$. The maximum is attained at $x = -\frac{(\mu^2 - 1)(e^{-t} - 1)}{(\mu^2 + 1)(e^{-t}-1)+2}$ if this value of $x$ satisfies $t_{+1}(x,\mu(1+x)) \leq -t \leq \hat t(x,\mu(1+x))$, and at $x = \frac{e^{-t}}{\mu\sqrt{e^{-2t} - 1}}-1$ if this value of $x$ satisfies $\hat t(x,\mu(1+x)) \leq -t \leq t^*(x,\mu(1+x))$. Straightforward calculation yields the maximizer claimed in the lemma.

The value of the maximum is obtained by evaluating the expression $B_{II}$ in the first case and the expression $B_{III}$ in the second case. Again straightforward calculation yields the value claimed in the lemma.
\end{proof}

\begin{figure}
\centering
\includegraphics[width=16.13cm,height=3.25cm]{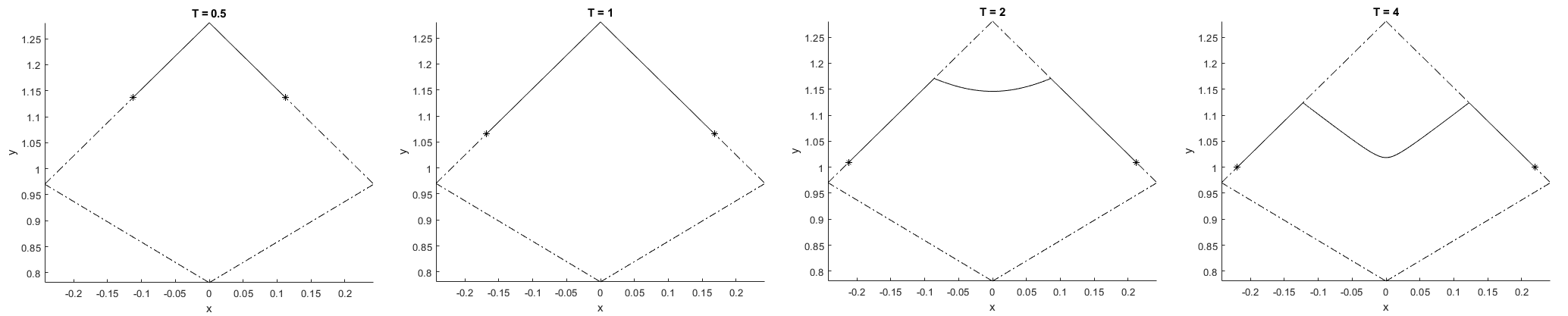}
\caption{Optimal trajectories of system \eqref{control_dynamics} with free end-points for different time horizons $T = t_f - t_i$. The optimal initial and terminal points are marked with stars. The parameter $\gamma$ equals $0.5$. The dash-dotted lines delimit the feasible region $X$.}
\label{fig:optimal_free}
\end{figure}

The optimal trajectory realizing the maximal value in Lemma \ref{lem:maximalB} is depicted in Fig.~\ref{fig:optimal_free}. For $t_f-t_i \leq \log\frac{\mu^2+1}{\mu^2-1}$ it consists of two arcs with control $u = \pm1$, respectively, and lies entirely on the boundary of the feasible set $X$ (top). For $t_f - t_i > \log\frac{\mu^2+1}{\mu^2-1}$ the optimal trajectory consists of three arcs with controls $+1,0,-1$, respectively. The first and third arc lie on the boundary of $X$, while the central arc is singular and crosses the interior of $X$ (bottom). The whole trajectory is symmetric about the vertical axis.

{\corollary \label{cor:upper_bound} Let $\Omega$ be a proper open segment of the projective line, and let $f$ be a convex non-degenerate centro-affine $C^2$ immersion of $\Omega$ into the interior of the cone $K_{\Omega} \subset \mathbb R^2$ over $\Omega$ which satisfies \eqref{Cbound} a.e.\ for some $\gamma > 0$ and is asymptotic to $\partial K_{\Omega}$. Then for every two points $a,b \in \Omega$ the Riemannian length $l^R(a,b)$ of the segment between $a,b$ in the centro-affine metric induced by $f$ is bounded above by \eqref{upper_bound_T}, where $T = d^H(a,b)$ is the Hilbert distance between $a,b$ in $\Omega$. }

\begin{proof}
Choosing a coordinate system in $\mathbb R^2$ such that $K_{\Omega} = \mathbb R_+^2$ and representing $f = e^{\alpha(t)}(e^t,e^{-t})$ by a $C^2$ function $\alpha(t)$ satisfying \eqref{Cbound3} a.e.\ reduces the problem of maximization of $l^R(a,b)$ to the optimal control problem \eqref{control_dynamics}--\eqref{state_constraints}. Application of Lemmas \ref{lem:h_upper_bound}--\ref{lem:maximalB} concludes the proof.
\end{proof}

In Corollary \ref{cor:upper_bound} we assumed that the immersion $f$ is asymptotic to the boundary of $K_{\Omega}$, while in Theorem \ref{thm:C_bound} this assumption is missing. In order to circumvent this difficulty we need the following lemmas.

{\lemma \label{lem:upper_bound_monotone} The difference $\delta = \max_{(x,y) \in X} B(-T,x,y) - T$ between \eqref{upper_bound_T} and $T$ is an increasing function of $T$ for $T > 0$. }

\begin{proof}
The derivative of the difference with respect to $T$ is given by
\[ \frac{d\delta}{dT} = \left\{ \begin{array}{rcl} \frac{(\mu - 1)(\mu + 1 - (\mu-1)e^T)}{\mu^2(e^T - 1) + e^T + 1},&\quad& e^T \leq \frac{\mu^2+1}{\mu^2-1}, \\
\frac{1-\sqrt{1-e^{-2T}}}{\sqrt{1-e^{-2T}}}, && e^T \geq \frac{\mu^2+1}{\mu^2-1}. \end{array} \right.
\]
For $e^T \leq \frac{\mu^2+1}{\mu^2-1}$ we have $\mu + 1 - (\mu-1)e^T \geq \mu + 1 - (\mu-1)\frac{\mu^2+1}{\mu^2-1} = \frac{2\mu}{\mu^2+1} > 0$, and for $T > 0$ we have $0 < \sqrt{1-e^{-2T}} < 1$. Hence the derivative is positive.
\end{proof}

{\lemma \label{lem:extension} Let $\alpha(t)$ be a function satisfying the conditions of Lemma \ref{lem:h_upper_bound} and such that $\ddot\alpha > \dot\alpha^2 - 1$. Then the centro-affine immersion $f$ defined by $\alpha$ is extendable to a convex non-degenerate centro-affine immersion $\tilde f$ of class $C^2$ which satisfies \eqref{Cbound} a.e.\ and is asymptotic to the boundary of some cone $\tilde K$ which contains $\mathbb R_+^2$. }

\begin{proof}
Let us consider the behaviour of $f(t) = e^{\alpha(t)}(e^t,e^{-t})$ as $t \to +\infty$. If $e^{\alpha(t)+t}$ is unbounded, then $f$ is asymptotic to the boundary ray $\rho$ of $\mathbb R_+^2$ generated by the vector $(1,0)$. Let us hence assume that $\lim_{t\to+\infty}(\alpha(t)+t) < +\infty$. Then $\liminf_{t \to \infty}\dot\alpha = -1$. However, by Lemma \ref{lem:h_upper_bound} we have $\ddot\alpha \leq (\mu^2+1)(\dot\alpha+1)^2 - 2(\dot\alpha+1)$ and hence $\ddot\alpha \leq 0$ whenever $\dot\alpha \leq -\frac{\mu^2-1}{\mu^2+1}$. It follows that $\lim_{t \to \infty}\dot\alpha = -1$.

Let us show that the immersion $f$ is transversal to the boundary ray $\rho$, i.e.,\ $\frac{\dot f_1}{\dot f_2} = e^{2t}\frac{\dot\alpha+1}{\dot\alpha-1}$ has a finite limit as $t \to +\infty$. It suffices to show that $\beta = e^{2t}(\dot\alpha+1)$ remains finite. Again by Lemma \ref{lem:h_upper_bound} we have
\[ \dot\beta \leq 2\beta + e^{2t}((\mu^2+1)(\dot\alpha+1)^2 - 2(\dot\alpha+1)) = (\mu^2+1)e^{-2t}\beta^2.
\]
It follows that with $\eta = \beta^{-1} - \frac{\mu^2+1}{2}e^{-2t}$ we have $\dot\eta = -\frac{\dot\beta}{\beta^2} + (\mu^2+1)e^{-2t} \geq 0$, and $\eta$ is increasing. However, $e^{2t}\eta(t) = \frac{1}{\dot\alpha+1} - \frac{\mu^2+1}{2}$ grows unbounded, and hence $\eta(t)$ eventually becomes positive. Thus $\beta = \frac{1}{\eta + \frac{\mu^2+1}{2}e^{-2t}}$ and consequently $\frac{\dot f_1}{\dot f_2}$ remain bounded. Then by convexity of the immersion the ratio $\frac{\dot f_1}{\dot f_2}$ must have a limit as $t \to +\infty$.

Let us show that the affine metric of $f$ has a limit. The coordinate $t$ becomes singular as the image $f(t)$ approaches the ray $\rho$. Let us therefore consider the non-singular coordinate $f_2 = e^{\alpha(t)-t}$. In this coordinate the affine metric is given by $(\ddot\alpha-\dot\alpha^2+1)\left(\frac{dt}{df_2}\right)^2 \leq \frac{\mu^2(1+\dot\alpha)^2}{((\dot\alpha-1)e^{\alpha(t)-t})^2} = \frac{\mu^2\beta^2}{((\dot\alpha-1)f_1)^2}$. The upper bound has a limit as $t \to +\infty$, and hence the metric remains bounded. By \eqref{Cbound} it is Lipschitz and has a well-defined limit.

On the other hand, this limit cannot be zero by condition \eqref{Cbound}. This can be seen from the equivalent form \eqref{Cbound3} (written down in a coordinate which is non-singular at the limit point), which ensures that $h$ cannot reach zero in finite time.

We may then extend $f$ continuously by a hyperbola branch which matches the limit values of the first two derivatives of $f$ at the intersection point with the ray $\rho$. This hyperbola will be asymptotic to a ray which defines the boundary of the cone $\tilde K$. The cubic form of an immersion defined by a quadric vanishes, and beyond $\rho$ the extension satisfies \eqref{Cbound} with $\gamma = 0$.

In the same way the extension of $f$ beyond the boundary ray of $\mathbb R_+^2$ generated by the vector $(0,1)$ is constructed, if $f$ is not already asymptotic to this ray.
\end{proof}

{\corollary \label{cor:1D_upper} Let $\Omega \subset \mathbb RP^1$ be a proper open segment of the projective line, and let $f: \Omega \to K_{\Omega}^o$ be a convex non-degenerate centro-affine lift of class $C^3$ into the interior of the cone over $\Omega$ which satisfies \eqref{Cbound}. For points $a,b \in \Omega$, let $T = d^H(a,b)$ be their Hilbert distance. Then the length $l^R(a,b)$ of the line segment between $a,b$ in the centro-affine metric defined by $f$ is strictly smaller than \eqref{upper_bound_T}, and this bound cannot be improved. }

\begin{proof}
By Lemma \ref{lem:extension} there exists a proper open segment $\tilde\Omega \supset \Omega$ of the projective line such that the immersion $f$ can be extended to a convex non-degenerate centro-affine lift $\tilde f$ of class $C^2$ of $\tilde\Omega$ into the interior of the cone $\tilde K$ over $\tilde\Omega$ which satisfies \eqref{Cbound} a.e.\ and is asymptotic to $\partial\tilde K$.

Let $\tilde d^H(a,b)$ be the Hilbert distance between the points $a,b$ with respect to $\tilde\Omega$. Then by Corollary \ref{cor:upper_bound} the Riemannian length $l^R(a,b)$ is upper bounded by \eqref{upper_bound_T} with $T = \tilde d^H(a,b)$. However, $\Omega \subset \tilde\Omega$ implies $\tilde d^H(a,b) \leq d^H(a,b)$. Since by Lemma \ref{lem:upper_bound_monotone} expression \eqref{upper_bound_T} is increasing with $T$, the length $l^R(a,b)$ is also upper bounded by \eqref{upper_bound_T} with $T = d^H(a,b)$. Moreover, the upper bound cannot be attained, because the optimal trajectory constructed in the proof of Lemma \ref{lem:maximalB} corresponds to an immersion which is only piece-wise $C^3$.

On the other hand, this optimal trajectory can be extended from the interval $[t_i,t_f]$ to $\mathbb R$ by applying control $u = 1$ for all $t < t_i$ and control $u = -1$ for all $t > t_f$. The extension then tends to the left-most point $(x,y) = \left(-\frac{\mu^2-1}{\mu^2+1},\frac{2\mu}{\mu^2+1}\right)$ of $X$ for $t \to -\infty$ and to the right-most point $(x,y) = \left(\frac{\mu^2-1}{\mu^2+1},\frac{2\mu}{\mu^2+1}\right)$ of $X$ for $t \to +\infty$. The corresponding centro-affine immersion into $\mathbb R_+^2$ is of class $C^2$ and piece-wise analytic, but can be approximated with arbitrary precision in the $C^2$ norm by $C^3$ immersions satisfying \eqref{Cbound}. Hence bound \eqref{upper_bound_T} cannot be improved.
\end{proof}

We may now return to Theorem \ref{thm:C_bound}. The first two inequalities are proven in a similar manner as for Theorem \ref{thm:C_libre}. Namely, the first one is just the inequality between the geodesic distance and the Riemannian length of a path linking $a,b$, while the second inequality is the upper bound on $l^R(a,b)$ obtained in Corollary \ref{cor:1D_upper} for the case $n = 1$ and which carries over to general dimension because the metric is centro-affine.

Let us prove the last inequality in Theorem \ref{thm:C_bound}. By Lemma \ref{lem:upper_bound_monotone} the difference $\delta$ mentioned in this lemma obeys
\begin{align*}
\delta(T) &< \lim_{s \to +\infty}\,\delta(s) = \lim_{s \to +\infty} \left( \log\frac{(\mu-1)(\sqrt{e^s+1} + \sqrt{e^s-1})}{(\mu+1)(\sqrt{e^s+1} - \sqrt{e^s-1})} + \frac{2\mu}{\mu^2+1}\log\frac{2\mu^2}{\mu^2-1} - s \right) \\ &= \log 2 - \log\frac{\mu + 1}{\mu - 1} + \frac{2\mu}{\mu^2+1}\log\frac{2\mu^2}{\mu^2-1}
\end{align*}
for every $T \geq 0$. Here in the second equality we used that $\frac{\sqrt{\kappa+1} + \sqrt{\kappa-1}}{\kappa(\sqrt{\kappa+1} - \sqrt{\kappa-1})} = 1 + \sqrt{1 - \frac{1}{\kappa^2}} \to 2$ as $\kappa \to +\infty$. Inserting $T = d^H(a,b)$ completes the proof of Theorem \ref{thm:C_bound}.

\section{Proof of Theorems \ref{thm:C_lower_bound} and \ref{thm:C_lower_bound_geodesic}} \label{sec:lower_bound}

In order to obtain a lower bound on the Riemannian length $l^R(a,b)$ we have to consider the optimal control problem
\begin{equation} \label{minimum_problem}
\begin{aligned}
\dot x = y^2 + x^2 - 1, &\quad \dot y = 2xy + u\gamma y^2, \qquad u \in [-1,1], \\
l^R(t_i,t_f) &= \int_{t_i}^{t_f} y(t)\,dt \to \inf, \\
\mu^{-1}(1+|x|) & \leq y \leq \mu(1-|x|),
\end{aligned}
\end{equation}
which is similar to \eqref{control_dynamics}--\eqref{state_constraints} with the difference that we now minimize \eqref{control_cost}. The proof is conducted along the same lines as that of Theorem \ref{thm:C_bound}, but the calculations turn out to be simpler because the optimal trajectories do not contain the singular arc and the optimal control is purely bang-bang (i.e., assuming only its extreme values). Assume the notations of the previous section.

{\lemma Assume notations \eqref{abcd}. Consider control problem \eqref{minimum_problem} with fixed initial point $(x(t_i),y(t_i)) = (x_i,y_i)$ and free terminal point $(x(t_f),y(t_f)) \in X$. Let $T = t_f - t_i$ be the time horizon. Then the optimal value of the problem is given by $B(-T,x_i,y_i)$, where the function $B: \mathbb R_- \times X \to \mathbb R_+$ is defined as follows.

Let $(x,y) \in X$ and $t \in \mathbb R_-$. If
\[ -t \leq t_{-1}(x,y) = \left\{ \begin{array}{rcl} \frac12\log\left( 1 + \frac{2\mu a_-}{c_-a_+} \right),&\quad& c_- > 0, \\ +\infty, && c_- = 0, \end{array} \right.
\]
then
\[ B(t,x,y) = \frac{\mu}{\mu^2+1}\log\frac{\mu(a_+e^{-2t} - b_-)}{c_-e^{-2t} + d_+}.
\]
If $-t > t_{-1}(x,y)$, then
\begin{align*}
B(t,x,y) &= \frac{\mu}{\mu^2+1}\log\left(\frac{a_+}{8\mu^3yc_-} \cdot \left[ (\mu^2 + 1)c_-a_+e^{-2t} + \right. \right. \\
&+ \left. \left. \mu^4y(1+x) + \mu(1+\mu^2)(y^2-x^2+1) - 6\mu^2y + y(1-x)\right] \vphantom{\frac12}\right).
\end{align*}
}

\begin{proof}
Let us show that $B(t,x,y)$ is the Bellman function of the problem. Denote the two expressions defining $B$ in the lemma by $B_I,B_{II}$.

\medskip

\emph{Consistency:} For $(x,y) \in X$ we clearly have $t_{-1}(x,y) \geq 0$.

\medskip

\emph{Continuity:} Inserting $t = -t_{-1}(x,y)$ into $B_I$ and $B_{II}$, we obtain the same expression
\[ \frac{\mu}{\mu^2+1} \log\frac{(\mu^2 - 1)a_+}{2\mu c_-}.
\]

\medskip

\emph{Initial value:} Inserting $t = 0$ into $B_I$, we obtain $B(0,x,y) \equiv 0$, which proves \eqref{Bellman_initial}.

\medskip

\emph{Bellman inequality:} Let us show \eqref{Bellman_equation}, i.e., that for every $u \in [-1,1]$ we have
\[
\frac{d}{dt}B(t,x,y) + y = \frac{\partial B}{\partial x} \cdot (y^2+x^2-1) + \frac{\partial B}{\partial y} \cdot (2xy+u\gamma y^2) + \frac{\partial B}{\partial t} + y \geq 0,
\]
with equality if $u = \hat u(t,x,y)$, where
\[ \hat u(t,x,y) = \left\{ \begin{array}{rcl} -1,&\quad& a_- > 0, \\ +1,&& a_- = 0. \end{array} \right.
\]

For $0 \leq -t < t_{-1}$ we have
\[ \frac{dB_I(t,x,y)}{dt} + y = \frac{y^2(\mu^2 - 1)(e^{-2t} - 1)(u + 1)((1 - x)e^{-2t} + 1 + x)}{(c_-e^{-2t} + d_+)(a_+(e^{-2t} - 1) + 2)} \geq 0,
\]
with equality if $u = -1$. Note that if $a_- = 0$, then $t_{-1} = 0$, and hence the inequality $-t < t_{-1}$ cannot hold.

For $-t \geq t_{-1}$ we have
\[ \frac{dB_{II}(t,x,y)}{dt} + y = \frac{y(\mu^2 - 1)(u + 1)a_-\left[a_+c_-^2(e^{-2t} - e^{2t_{-1}}) + 2\mu y((\mu^2 - 1)(1 - x) + \mu c_-)\right]}{c_-a_+\left((\mu^2 + 1)c_-a_+(e^{-2t} - e^{2t_{-1}}) + 4\mu^2y(\mu^2 - 1)\right)} \geq 0,
\]
with equality if $u = -1$ or $a_- = 0$.

\medskip

\emph{Feasibility:} The only boundary segments through which a trajectory can escape $X$ are the upper right and the lower left one. On the lower left segment we have $a_- = 0$ and hence the optimal control is $\hat u = +1$. The trajectory then moves along the boundary segment. On the upper right segment we have  $t_{-1} = +\infty$ and $a_- > 0$. Therefore the optimal control is $\hat u = -1$ and the trajectory again moves along the boundary segment.

\medskip

Thus $B$ is indeed the Bellman function, and the optimal value is achieved by applying the control $u = \hat u$.
%
\end{proof}

\begin{figure}
\centering
\includegraphics[width=10.68cm,height=4.17cm]{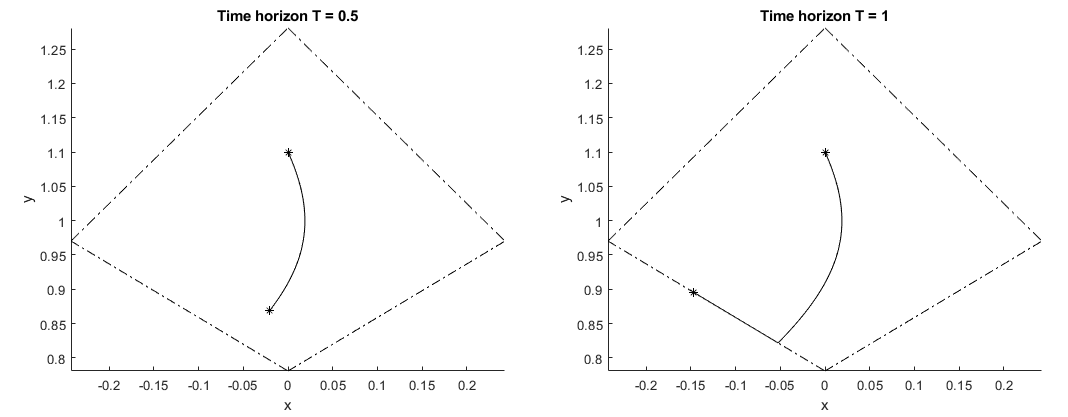}
\caption{Optimal trajectories of system \eqref{minimum_problem} with initial point $(x,y) = (0,1.1)$ for different time horizons $T = t_f - t_i$. The initial and terminal points are marked with stars. The parameter $\gamma$ equals $0.5$. The dash-dotted lines delimit the feasible region given by the state constraints in \eqref{minimum_problem}.}
\label{fig:optimal_synthesis_min}
\end{figure}

The optimal solutions obtained by application of control $\hat u$ are structured as follows.

If $T \leq t_{-1}(x_i,y_i)$, where $t_{-1}$ is the time needed to reach the lower left boundary segment of the feasible set with control $u \equiv -1$, then the optimal control is given by $\hat u \equiv -1$ on the whole trajectory (see Fig.~\ref{fig:optimal_synthesis_min}, left).

If $T > t_{-1}(x_i,y_i)$, then on the optimal trajectory the control $u \equiv -1$ is optimal for $t \in [t_i,t_i+t_{-1}(x_i,y_i))$, that is up to the point when the trajectory reaches the boundary of the feasible set. For $t \in (t_i+t_{-1}(x_i,y_i),t_f]$ the control $u \equiv +1$ is optimal and the trajectory moves along the boundary of the feasible set (see Fig.~\ref{fig:optimal_synthesis_min}, right).

\medskip

In order to find the optimal value of problem \eqref{minimum_problem} with free initial and terminal points, we have to minimize the Bellman function $B(-T,x,y)$ over $(x,y) \in X$ for fixed $T$.

{\lemma \label{lem:minimalB} The minimum $\min_{(x,y) \in X}B(-T,x,y)$ is attained at
\[ (x,y) = \left(\frac{(e^T - 1)(\mu^2 - 1)}{e^T(\mu^2 + 1) + \mu^2 - 1},\frac{2\mu e^T}{e^T(\mu^2 + 1) + \mu^2 - 1}\right).
\]
The corresponding value of the minimum is given by
\begin{equation} \label{lower_bound}
\min_{(x,y) \in X}B(-T,x,y) = \frac{2\mu}{\mu^2 + 1}\log\frac{e^T(\mu^2 + 1) + \mu^2 - 1}{2\mu^2}.
\end{equation}
}

\begin{proof}
We shall again parameterize $X$ by the variables $w = \frac{y}{1-x}$, $z = \frac{y}{1+x}$. Set $t = -T$.

Let us show that $\frac{\partial B}{\partial z} \geq 0$. For $-t \leq t_{-1}(x,y)$ we have
\[ \frac{\partial B_I}{\partial z} = \frac{\mu(x + 1)^2(e^{-2t} - 1)}{(c_-e^{-2t} + d_+)(a_+(e^{-2t} - 1) + 2)} \geq 0.
\]
For $-t \geq t_{-1}(x,y)$ we have
\[ \frac{\partial B_{II}}{\partial z} = \frac{\mu(e^{-2t} - 1)(x + 1)^2c_-a_+}{2y\left((\mu^2 + 1)c_-a_+(e^{-2t} - e^{2t_{-1}}) + 4\mu^2y(\mu^2 - 1)\right)} \geq 0.
\]
Hence the minimum of $B$ is achieved at $z = \mu^{-1}$. This corresponds to the lower right boundary segment of $X$.

We now compute the derivative $\frac{\partial B}{\partial w}$ on this segment. On this segment $x \in \left(0,\frac{\mu^2-1}{\mu^2+1}\right)$ and $e^{2t_{-1}} = \frac{(\mu^2 - 1)(x + 1)}{\beta}$, where $\beta = \mu^2 - 1 - x(\mu^2 + 1) > 0$.

For $-t < t_{-1}(x,y)$ we have
\[ \frac{\partial B_I}{\partial w} = \frac{\mu(e^{-2t} - 1)}{\beta e^{-2t} + (\mu^2 + 1)(x + 1)} \geq 0.
\]
At $x = 0$ we have $t_{-1} = 0$, and hence the minimum cannot be attained for $-t < t_{-1}(x,y)$.

For $-t \geq t_{-1}(x,y)$ we have
\[ \frac{\partial B_{II}}{\partial w} = \frac{\mu\left[(\mu^2 - 1)^2(x + 1)^2 - \beta^2e^{-2t}\right]}{(x + 1)\beta\left[(\mu^2 + 1)\beta(e^{-2t} - e^{2t_{-1}}) + 2\mu^2(\mu^2 - 1)(x + 1)\right]}.
\]
Hence $\frac{\partial B_{II}}{\partial w} \geq 0$ if $x \geq \frac{(e^{-t} - 1)(\mu^2 - 1)}{e^{-t}(\mu^2 + 1) + \mu^2 - 1}$ and $\frac{\partial B_{II}}{\partial w} \leq 0$ if $x \leq \frac{(e^{-t} - 1)(\mu^2 - 1)}{e^{-t}(\mu^2 + 1) + \mu^2 - 1}$.

It follows that the minimum is attained at
\[ (x,y) = \left( \frac{(e^{-t} - 1)(\mu^2 - 1)}{e^{-t}(\mu^2 + 1) + \mu^2 - 1},\mu^{-1}(1+x) \right).
\]
The value of the minimum is obtained by evaluating the expression $B_{II}$ at this point. 
\end{proof}

\begin{figure}
\centering
\includegraphics[width=10.55cm,height=4.18cm]{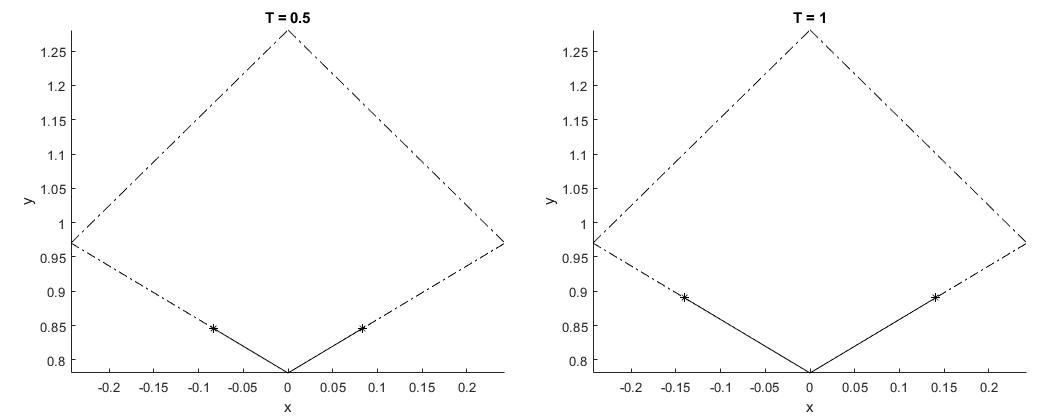}
\caption{Optimal trajectories of system \eqref{minimum_problem} with free end-points for different time horizons $T = t_f - t_i$. The optimal initial and terminal points are marked with stars. The parameter $\gamma$ equals $0.5$. The dash-dotted lines delimit the feasible region given by the state constraints in \eqref{minimum_problem}.}
\label{fig:optimal_free_min}
\end{figure}

The optimal trajectory realizing the minimal value in Lemma \ref{lem:minimalB} is depicted in Fig.~\ref{fig:optimal_free_min}. It consists of two arcs with control $u = \mp1$, respectively, and lies entirely on the boundary of the feasible set $X$. It is symmetric about the vertical axis.

The solutions can be extended from the time interval $[t_i,t_f]$ to $\mathbb R$ by applying control $u = -1$ for all $t < t_i$ and control $u = +1$ for all $t > t_f$. The corresponding trajectory then tends to the right-most point $(x,y) = \left(\frac{\mu^2-1}{\mu^2+1},\frac{2\mu}{\mu^2+1}\right)$ of $X$ for $t \to -\infty$ and to the left-most point $(x,y) = \left(-\frac{\mu^2-1}{\mu^2+1},\frac{2\mu}{\mu^2+1}\right)$ of $X$ for $t \to +\infty$. The corresponding centro-affine immersion into $\mathbb R_+^2$ is of class $C^2$ and piece-wise analytic, but can be approximated with arbitrary precision in the $C^2$ norm by $C^3$ immersions satisfying \eqref{Cbound}. Hence the lower bound \eqref{lower_bound} on $l^R(-T,0)$ cannot be attained by $C^3$ immersions, but is nevertheless sharp.

\medskip

Let us now prove Theorem \ref{thm:C_lower_bound}. The first inequality is just the bound \eqref{lower_bound}, which carries over to the case of general dimension as in the proof of Theorem \ref{thm:C_libre}. The second inequality comes from the relation
\[ \log\frac{e^{\kappa}(\mu^2+1) + \mu^2-1}{2\mu^2} - \left( \kappa - \log\frac{2\mu^2}{\mu^2+1} \right) = \log\left(1 + e^{-\kappa}\frac{\mu^2-1}{\mu^2+1}\right) > 0,
\]
which also becomes sharp as $\kappa \to +\infty$.

\medskip

Let us now prove Theorem \ref{thm:C_lower_bound_geodesic}. Let $a,b \in \Omega$ be arbitrary points, and let $\sigma$ be the Riemannian geodesic linking these points. The Riemannian length of $\sigma$ is by definition equal to $d^R(a,b)$. Let $l$ be the length of the curve $\sigma$ in the Hilbert metric. Then $d^H(a,b) \leq l$, because straight lines are the shortest paths in the Hilbert metric. Summing the first inequality in Corollary \ref{cor:infinitesimal_length} over increasingly finer partitions of the curve $\sigma$ we obtain $\mu^{-1}l \leq d^R(a,b)$ in the limit. Hence $\mu^{-1}d^H(a,b) \leq d^R(a,b)$, yielding the first inequality in Theorem \ref{thm:C_lower_bound_geodesic}.

On the other hand, combining the second inequality in Corollary \ref{cor:infinitesimal_length} with the relation $d^R(a,b) \leq l^R(a,b)$ we obtain $d^R(a,b) \leq \mu d^H(a,b)$, which is the second inequality in Theorem \ref{thm:C_lower_bound_geodesic}. This completes the proof of Theorem \ref{thm:C_lower_bound_geodesic}.

\bibliography{affine_geometry,opt_control,hilbertGeometry}
\bibliographystyle{plain}

\end{document}